\documentclass[12pt,leqno]{article} 
\usepackage{amsmath,amssymb,latexsym,theorem,epsfig}
\usepackage[all]{xy}
\usepackage{psfrag}

\makeatletter
\let\@internalcite\cite
\def\cite{\def\citeauthoryear##1##2{##1{##2}}\@internalcite}
\def\@biblabel#1{\def\citeauthoryear##1##2{##1{##2}}[#1]\hfill}
\makeatother

\theoremstyle{plain} 

\newtheorem{thm}{Theorem}
\newtheorem{theo}{Theorem}[section]
\newtheorem{lem}[theo]{Lemma}
\newtheorem{cor}[theo]{Corollary}
\newtheorem{con}[theo]{Conjecture}
\newtheorem{defi}[theo]{Definition}

\newtheorem{prop-defi}[theo]{Proposition-Definition}
\newtheorem{lemma-defi}[theo]{Lemma-Definition}
{\theorembodyfont{\rmfamily} \newtheorem{rem}[theo]{Remark}}
{\theorembodyfont{\rmfamily} }

\newcommand{\op}[1]{\operatorname{#1}}
\newcommand{\eop}{\ \hfill $\Box$}
\newcommand{\balpha}{\boldsymbol{\alpha}}
\newcommand{\ba}{\boldsymbol{a}}
\newcommand{\be}{\boldsymbol{e}}
\newcommand{\bsigma}{\boldsymbol{\sigma}}
\newcommand{\bmu}{\boldsymbol{\mu}}
\newcommand{\bB}{\boldsymbol{B}}
\newcommand{\bC}{\boldsymbol{C}}
\newcommand{\bT}{\boldsymbol{T}}

\newcommand{\bone}{\boldsymbol{{\mathit 1}}}
\newcommand{\liesp}{{\mathfrak s}{\mathfrak p}}
\newcommand{\lieso}{{\mathfrak s}{\mathfrak o}}
\newcommand{\liesl}{{\mathfrak s}{\mathfrak l}}
\newcommand{\bbE}{{\mathbb E}}
\newcommand{\bbD}{{\mathbb D}}
\newcommand{\bbG}{{\mathbb G}}

\newcommand{\monb}[1]{\op{mon}_{B}^{#1}}
\newcommand{\mondr}[1]{\op{mon}_{DR}^{#1}}

\numberwithin{equation}{section}

\begin{document}
\centerline{\LARGE\sf Density of monodromy actions on non-abelian
cohomology}

\

\bigskip
\bigskip

\noindent
{\large Ludmil Katzarkov}\footnote{Partially supported by NSF Career
Award DMS-9875383 and A.P. Sloan research fellowship} \\
University of California at Irvine \\ 
Irvine, CA 92697, USA

\medskip

\noindent
{\large Tony Pantev}\footnote{Partially supported by 
NSF Grant DMS-9800790 and A.P. Sloan research fellowship} \\
University of  Pennsylvania \\ 
209 South 33rd street \\ 
Philadelphia, PA 19104-6395, USA

\medskip

\noindent
{\large  Carlos Simpson}  \\
CRNS - Universit\'{e} de Nice-Sophia Antipolis \\
Parc Valrose, 06108 Nice Cedex 2, FRANCE

\

\bigskip

\

\begin{abstract} In this paper we study 
the monodromy action on the first Betti and de Rham 
non-abelian cohomology 
arising from a family of smooth curves. We describe sufficient
conditions for the existence of a Zariski dense monodromy orbit. In
particular we show that for a Lefschetz pencil of sufficiently high
degree the monodromy action is dense.
\end{abstract}

\tableofcontents

\section{Introduction} \label{sec-introduction}

We work in the category of quasi-projective schemes over ${\mathbb C}$.
Let $f: X \to B$ be a smooth projective morphism with connected 
fibers of dimension one and genus at least two. Fix a base point $o \in B$ 
and let $X_{o}$ be the corresponding fiber of $f$. In this paper we
study the monodromy action of $\pi_{1}(B,o)$ on the degree one non-abelian
Betti and de Rham cohomology of $X_{o}$. Let $\pi_{1}(X_{o})$ denote
the abstract fundamental group of $X_{o}$ and let 
\[
\op{mon} : \pi_{1}(B,o) \to \op{Out}(\pi_{1}(X_{o}))
\]
be the geometric monodromy representation of the family $f : X \to B$.
For any positive integer $n$ let
\[
\xymatrix@R=4pt{
\monb{n} : \pi_{1}(B,o) \ar[r] &
\op{Aut}(H^{1}_{B}(X_{o},\op{GL}(n,{\mathbb C})))}  
\]
be the induced monodromy action on the non-abelian Betti
cohomology with coefficients in $\op{GL}(n,{\mathbb C})$. 

The de Rham object which corresponds to $\monb{n}$ is the Gauss-Manin
connection $\nabla _{DR}^n$ on the relative de Rham cohomology stack
$H^1 _{DR}(X/B, \op{GL}(n,{\mathbb C}))$.  While the non-abelian Betti
and de Rham cohomology are most naturally viewed as stacks, for the
purposes of the present paper it suffices to work with the
corresponding coarse moduli spaces. To indicate that we will write
$M_{B}(X_{o},n)$ and $M_{DR}(X/B,n)$ rather than
$H^{1}_{B}(X_{o},\op{GL}(n,{\mathbb C}))$ and
$H^{1}_{DR}(X/B,\op{GL}(n,{\mathbb C}))$. Concretely $M_{B}(X_{o},n)$
denotes the moduli space of (semisimplifications of) represenations of
$\pi_{1}(X_{o})$ in $\op{GL}(n,{\mathbb C})$ and $M_{DR}(X_{o},n)$
denotes the moduli space of rank $n$ algebraic local systems of
complex vector spaces on $X_{o}$. The total space
$M_{DR}(X/B,n)\rightarrow B$ is a quasiprojective variety over $B$
whose fiber over the point $o$ is $M_{DR}(X_{o},n)$.

For a loop $\gamma \in
\pi_{1}(B,o)$ the action of $\monb{n}(\gamma)$ on $M_{B}(X_{o},n)$ is
given by composing a $n$-dimensional represenation of
$\pi_{1}(X_{o})$ with some lift of the outer automorphism
$\op{mon}(\gamma)$ of $\pi_{1}(X_{o})$. This gives a well defined
action on conjugacy classes of representations of
$\pi_{1}(X_{o})$ and results in an algebraic automorphism
$\monb{n}(\gamma) : M_{B}(X_{o},n) \to M_{B}(X_{o},n)$. 

There is an analytic
action $\mondr{n}$ of $\pi _1(B,o)$ on $M_{DR}(X_{o},n)$ 
which is most naturally described
through the Riemann-Hilbert correspondence (see
e.g. \cite{deligne-ode-book},
\cite[Section~7]{simpson-moduli2}). It
is shown in  \cite[Section~7]{simpson-moduli2} that the passage from
a local system to its monodromy representation induces an 
isomorphism of analytic spaces 
\[
\tau : M_{DR}(X_{o},n)^{\op{an}}
\widetilde{\to} M_{B}(X_{o},n)^{\op{an}}.
\] 
Now given $\gamma$ we can
define an analytic automorphism of $\mondr{n}(\gamma)$ of
$M_{DR}(X_{o},n)$ by putting $\mondr{n}(\gamma) = \tau^{-1}\circ
\monb{n}(\gamma) \circ \tau$. This analytic action is the monodromy of 
the algebraic nonabelian Gauss-Manin connection $\nabla ^n_{DR}$ on the 
total space
$M_{DR}(X/B,n)$ \cite[Section~8]{simpson-moduli2}.

It is natural to try to understand the complexity of the  algebraic
(respectively analytic) action of $\pi_{1}(B,o)$ on $M_{B}(X_{o},n)$
(respectively  
$M_{DR}(X_{o},n)$) by measuring in some way the size of the 
$\pi_{1}(B,o)$-orbits on $M_{B}(X_{o},n)$ and
$M_{DR}(X_{o},n)$. Analogous questions concerning  the monodromy
action of $\pi_{1}(B,o)$ on spaces of special representations of
$\pi_{1}(X_{o})$ (e.g. real representations, representations with
compact image, projective structures, etc.) have been the focus of
active research in the recent years \cite{goldman-ergodic},
\cite{mcmullen}, \cite{gkm}. In that direction the result most relevant to
our setup is  a theorem of W.~Goldman who showed in
\cite{goldman-ergodic} that the mapping class group acts ergodically
on the space of all representations of $\pi_{1}(X_{o})$ into
$SU(2)$. Unfortunately, Goldman's proof does not generalize to the
case of representations into $SU(n)$ for $n > 2$ and we do not know whether
the mapping class group of $X_{o}$ still acts ergodically on the space
of such representations. Instead of pursuing the ergodicity question
in its full generality we chose to work with non-abelian cohomology
with complex coefficients. This allows us to use the algebraic
(respectively analytic) nature of the monodromy action on
$M_{B}(X_{o},n)$ (respectively $M_{DR}(X_{o},n)$) and to describe the
size of the monodromy orbits on those spaces in geometric, rather than
measure-theoretic terms.

Our first result is of essentially topological nature. Before we state
it we will need to introduce some notation. Let as before $f : X \to B$ be a
smooth holomorphic family of genus $g$ curves and let $o \in B$ be a
base point. Let $X_{o}$ be the fiber of $f$ over $o$ and let
\[
\op{map} : \pi_{1}(B,o) \to \op{Map}(X_{o}) \subset \op{Out}(\pi_{1}(X_{o})
\]
be the corresponding geometric monodromy representation.
Here $\op{Map}(X_{o})$ denotes the mapping class group of $X_{o}$. By
definition $\op{Map}(X_{o}) := \pi_{0}(\op{Diff}^{+}(X_{o}))$ is the
group of connected components of the group of orientation preserving
diffeomorphisms of $X_{o}$. Alternatively $\op{Map}(X_{o})$ can be
identified with the subgroup of index two in $\op{Out}(\pi_{1}(X_{o})$
consisting of all outer automorphisms which act trivially on
$H^{2}(X_{o},{\mathbb Z})$. Fix some topological double cover $\nu :
X_{o} \to {\mathbb P}^{1}$ and let $\imath \in \op{Map}(X_{o})$ be the
mapping class of the covering involution. The {\em hyperelliptic
mapping class group} of $X_{o}$ is defined to be the centralizer
$\Delta(X_{o})$ of $\imath$ in $\op{Map}(X_{o})$:
\[
\Delta(X_{o}) := \{ \phi \in \op{Map}(X_{o}) | \phi\imath \phi^{-1} =
\imath \}.
\]
Note that the definiton of $\Delta(X_{o})$ depends on the choice of
the double cover $\nu$ and so the hyperelliptic mapping class group is
defined as a subgroup in $\op{Map}(X_{o})$ only up to conjugation.
We will say that
the geometric monodromy of $f$ {\em dominates the hyperelliptic monodromy} if
up to  conjugation in
$\op{Map}(X_{o})$ the monodromy group $\op{mon}(\pi_{1}(B,o)) \subset 
\op{Map}(X_{o})$ contains a subgroup of finite
index in $\Delta(X_{o})$. 

\begin{thm} \label{thm-finite-index}
Assume that the monodromy of $f$ dominates the hyperelliptic
monodromy, e.g. assume that the image $\op{mon}(\pi_{1}(B,o)) \subset
\op{Out}(\pi_{1}(X_{o}))$ is of finite index in
$\op{Out}(\pi_{1}(X_{o}))$. Then there exists a positive integer
$g_{0}$ so that if $g \geq g_{0}$ and $n$ is any fixed odd integer, we
have:
\begin{description}
\item[{\em (i)}] There is no meromorphic function on
$M_{B}(X_{o},n)^{\op{an}}$ which is invariant under 
the action of $\monb{n}(\pi_{1}(B,o))$ 
(equivalently there is no meromorphic function on 
$M_{DR}(X/B,n)^{\op{an}}$ which is 
$\nabla ^n_{DR}$-invariant);
\item[{\em (ii)}] In the case of $M_{B}(X_{o},n)$, considered with its
natural structure of an affine algebraic variety,  there exists a point $x_{B}
\in M_{B}(X_{o},n)$ so that the orbit
\[
\monb{n}(\pi_{1}(B,o))\cdot x_{B}  \subset M_{B}(X_{o},n)
\] 
is Zariski dense in $M_{B}(X_{o},n)$;
or  in the case of 
$M_{DR}(X/B, n)\rightarrow B$ there is
a leaf of the foliation defined by $\nabla ^n_{DR}$ which is Zariski dense
in the algebraic Zariski topology.
\end{description}
\end{thm}

This theorem suggests that for families $f : X\to B$ with a ``large
enough'' geometric monodromy one should expect Zariski dense
monodromy actions on non-abelian cohomology. Geometrically families
with large monodromy naturally arise from hyperplane sections and
Lefschetz fibrations. In this context we prove the following

\begin{thm} \label{thm-lefschetz-pencil}
Let $Z$ be a smooth projective surface with
$b_{1}(Z)=0$. Let ${\mathcal O}_{Z}(1)$ be an ample line bundle on $Z$
and let $n > 1$ be a fixed odd
integer. Then there exists a positive integer $\ell$ (depending only
on $Z$ and ${\mathcal O}_{Z}(1)$), such that for every $k \geq \ell$ 
and for every Lefschetz fibration $f: \widehat{Z} 
\to {\mathbb P}^{1}$ in the linear system $|{\mathcal O}_{Z}(k)|$ we
have:
\begin{description}
\item[{\em (i)}] There is no meromorphic function on
$M_{B}(Z_{o},n)^{\op{an}}$ which is invariant under 
the action of $\monb{n}(\pi_{1}({\mathbb P}^{1}\setminus \{ p_{1}, \ldots, 
p_{\mu} \},o))$ 
(equivalently there is no meromorphic function on 
$M_{DR}(\widehat{Z}/ {\mathbb P}^{1}\setminus \{ p_{1}, \ldots, 
p_{\mu} \},n)^{\op{an}}$ which is 
$\nabla^n_{DR}$-invariant);
\item[{\em (ii)}] In the case of $M_{B}(Z_{o},n)$,  there exist a point $x_{B}
\in M_{B}(Z_{o},n)$ so that the orbit
\[
\monb{n}(\pi_{1}({\mathbb P}^{1}\setminus \{ p_{1}, \ldots, 
p_{\mu} \},o))\cdot x_{B}  \subset M_{B}(Z_{o},n)
\] 
is Zariski dense in $M_{B}(Z_{o},n)$;
or in the case of the space
$M_{DR}(\widehat{Z}/{\mathbb P}^{1}\setminus \{ p_{1}, \ldots, 
p_{\mu} \}, n)$ the foliation defined by $\nabla ^n_{DR}$ has a Zariski dense
leaf.
\end{description}
Here as
usual $\widehat{Z}$ is the blow-up of $Z$ at the base points of the pencil and
$p_{1}, \ldots p_{\mu} \in {\mathbb P}^{1}$ are the points where the
map $\widehat{Z} \to {\mathbb P}^{1}$ is not submersive.
\end{thm}

These statements can be viewed as nonabelian analogues of Deligne's 
irreducibility theorem \cite[Section~4.4]{deligne-weil2} and
\cite{janssen}, which
asserts that the monodromy group on the first cohomology of a
Lefschetz pencil of curves is a subgroup of finite index in the full 
symplectic group of the lattice of vanishing cycles.

\

\medskip

\noindent
The paper is organized as follows. In Section~\ref{sec-open-orbits} 
we examine the action of a finitely generated group on an affine
algebraic variety. We show how the existence of a Zariski dense orbit
can be deduced from the existence of an open orbit for 
the linearized action on
the tangent space at a fixed point. Section~\ref{sec-Schrodinger} describes 
a particular
point in the moduli space of representations of the fundamental group
of a curve $X_{o}$, which corresponds to the
Schr\"{o}dinger representation of a suitably chosen finite dihedral 
Heisenberg group. This
point is smooth and fixed by a subgroup of finite index in the
monodromy.  Moreover the tangent space of the moduli of representations at the
`Schr\"{o}dinger point' is naturally identified with the first cohomology of an
etale cover $Y_{o}$ of $X_{o}$.  Finally in
Section~\ref{sec-irreducible} we discuss the necessity of the
hypotheses of Theorem~\ref{thm-finite-index} and
\ref{thm-lefschetz-pencil} for the existence of a dense monodromy
orbit. We conjecture that the density holds under very mild
assumptions and give some additional evidence supporting the
conjecture.

\bigskip

\noindent
{\bf Acknowledgments.} We would like to thank A.~Beilinson, P.~Deligne
and I.~Smith for some very enlightening conversations. We are very
grateful to B.~Toen, C.~Walter and to the referee for pointing out
some mistakes in a preliminary version of this paper and for
suggesting improvements in the exposition. The second and third
authors would like to thank UC Irvine for their hospitality during the
month of July of 1998, when most of the ideas for the present paper
took shape. The first and the second author wish to thank the RiP
program of the Mathematical Forschunginstitut Oberwolfach and the
Volkswagen-Stiftung for their support and the excellent working
conditions during two weeks in the Summer of 1999, when a substantial
part of this work was done.

\section{Preliminary reductions} \label{sec-preliminaries}

We start with some general results about linear group actions on
algebraic varieties, which will allow us to localize at a point 
the Zariski density property of an action.

\subsection{Open orbits and dense actions} \label{sec-open-orbits}

Suppose $M$ is a reduced irreducible affine scheme of finite type over
${\mathbb C}$.  Write $M=\op{Spec}(A)$ and let $\Gamma$ be a finitely
presented discrete group acting on $M$ by algebraic automorphisms.
Thus $\Gamma$ acts on $A$ by ${\mathbb C}$-algebra automorphisms.

\begin{lem} \label{lem-generic-point}
If there is one point in $M$ whose orbit under $\Gamma$ is Zariski dense,
then there is a countable union of proper closed subvarieties of $M$ such
that for any $x$ not on this countable union, the orbit of $x$ is 
Zariski dense. 
\end{lem}
{\bf Proof.}  The fact that $A$ is of finite type over ${\mathbb C}$
means that we have a surjection from a polynomial ring to $A$. In
particular (doing an enumeration of the monomial basis for this
polynomial ring) we obtain a filtration of $A$: 
\[
A = \bigcup_{i=0}^{\infty} A_i
\]
by finite-dimensional
sub-${\mathbb C}$-vector spaces $A_i \subset A$.

For each pair of integers $i,k\geq 0$ let $G_{i,k}'$ denote the
space of $k$-tuples of elements of $A_i$. It is a finite dimensional affine
space. For $V\in G_{i,k}'$ let $I_A(V)$ denote the ideal in $A$ generated
by $V=(v_1,\ldots , v_k)$ and let $Z_V\subset M$ denote the reduced 
closed subvariety of $M$
defined by $V$. 

There is a closed algebraic subset
\[
{\mathcal Z}_{i,k}' \subset G_{i,k}' \times M
\] 
such that for each $V\in G_{i,k}'$ the fiber over $V$ is equal to
$Z_V$.  To see this, note that the coordinate ring of $G_{i,k}'$ is
the symmetric algebra on the dual of $A_i ^{\oplus k}$.  The
coordinate ring of $G_{i,k}'\times M$ is thus the tensor product and
the projections-inclusions $A_i^{\oplus k} \rightarrow A$ can be
viewed as elements of this coordinate ring; they generate the ideal of
the closed subset ${\mathcal Z}_{i,k}'$.  Let $G_{i,k} \subset
G'_{i,k}$ be the complement of the origin and let ${\mathcal Z}_{i,k}$
be the inverse image of $G_{i,k}$. The family
\[
{\mathcal Z}_{i,k} \subset G_{i,k} \times M
\] 
parameterizes precisely the closed proper subvarieties of $M$ which are
cut out by ideals generated by $k$-tuples of elements in
$A_{i}$. (Note that since $\op{Spec}(A)$ is reduced and irreducible,
any non zero $k$-tuple generates the ideal of a proper subvariety).

For each $\gamma \in \Gamma$ we obtain the translate
\[
\gamma {\mathcal Z}_{i,k} \subset G_{i,k} \times M.
\]
Let $G_{i,k} (\gamma ) \subset G_{i,k}$ denote the subset of points $V$ such 
that $\gamma Z_V = Z_V$.  We claim that this is a constructible
subset. It may be described as the set of points $V\in G_{i,k}$ such that
the intersection 
$\gamma {\cal Z}_{i,k}\cap {\cal Z}_{i,k}\cap (\{ V\} \times M)$
contains ${\cal Z}_{i,k}\cap (\{ V\} \times M)$.  This set is the complement
in $G_{i,k}$ of the  image of the map 
\[
({\mathcal Z}_{i,k}-(\gamma {\cal Z}_{i,k}\cap {\cal
Z}_{i,k}))\rightarrow G_{i,k}
\]
so it is constructible.

Now as $\gamma_j$ runs through a finite set of generators we obtain
the intersection of this finite collection of subsets, which is again
a constructible subset 
\[
G_{i,k}(\Gamma ):= \bigcap _{j}G_{i,k}(\gamma _j)\subset G_{i,k}.
\]
Each $V$ in $G_{i,k}(\Gamma )$ 
defines a $\Gamma$-invariant closed proper subvariety
$Z_{V} \subset M$, and conversely
it is clear that any $\Gamma$-invariant closed proper subvariety of $M$ 
appears as a $Z_V$ for some $i,k$ and some $V\in G_{i,k}(\Gamma )$.

Let $N_{i,k}$ be the union of the points contained in all of the
the subvarieties $Z_V$ for all $V\in G_{i,k}(\Gamma )$. This is again
a constructible set since it is the image of the
projection 
\[
{\mathcal Z}_{i,k}\times _{G_{i,k}}G_{i,k}(\Gamma ) \rightarrow M.
\]
If a point $x\in M$ is contained in any proper closed
$\Gamma$-invariant subvariety then it is contained in some
$N_{i,k}$. Note also that $N_{i,k}$ is $\Gamma$-invariant; indeed it
is a union of the $\Gamma$-invariant subsets $Z_V$ for all the $V\in
G_{i,k}(\Gamma )$.

We now have two possibilities: either
\begin{description}
\item[(a)] one of the constructible subsets $N_{i,k}$ is dense in $M$; or 
\item[(b)] all of the constructible subsets $N_{i,k}$ 
are contained in proper closed subvarieties $\overline{N}_{i,k}$.
\end{description}
In case (b), we obtain a countable union of closed subvarieties
$\bigcup _{i,k}\overline{N}_{i,k}$ such that if $x\in M$ is a point
which is not in this countable union, then $x$ is never contained in a
proper closed $\Gamma$-invariant subvariety.

In case (a) we claim that no point has a Zariski dense orbit.  Indeed 
the complement of the constructible set which is dense, has a closure which
is itself a proper Zariski closed and $\Gamma$-invariant subvariety 
(note that all of our constructible sets were $\Gamma$-invariant).
Thus any point here has non Zariski dense orbit.  On the other hand,
any point in the complement of this closed set is in the open interior
of the constructible set in question, so it is by definition in the image of
one of the $Z_V$, i.e. it is in a proper $\Gamma$-invariant closed subvariety.
This proves that in case (a) no point has a Zariski dense orbit.

Assume now that there is some point in $M$ which has a Zariski dense
orbit, then we are in case (b), so there is a countable union of
closed subvarieties such that if $x$ is not in here then $x$ has
Zariski dense orbit. The lemma is proven. \eop

\bigskip

\noindent
We are now ready to prove the main result of this section

\begin{theo} \label{theo-alternative}
Suppose that $\Gamma$ acts on an irreducible complex affine algebraic
variety $M$. There are two possibilities: either
\begin{description}
\item[(1)] there exists a nonconstant $\Gamma$-invariant meromorphic
function; or
\item[(2)] there exists a point $x\in M$ with $\Gamma x$ Zariski-dense in $M$.
\end{description}
\end{theo}
{\bf Proof.} Assume that no point in $M$ has a Zariski dense
orbit. This means that we are in the situation (a) discussed in  the
proof of Lemma~\ref{lem-generic-point}. In other words, there exist
integers $i, k \geq 0$ so that the $\Gamma$-invariant constructible
set $N_{i,k} \subset M$ is dense in $M$. To simplify notation put $S :=
G_{i,k}(\Gamma)$ and ${\mathcal Z} := {\mathcal
Z}_{i,k}\times_{G_{i,k}} G_{i,k}(\Gamma)$. By construction $S$ and
${\mathcal Z}$ are both schemes of finite type over ${\mathbb C}$ and the
natural maps
\[
\xymatrix{
{\mathcal Z} \ar[d] \ar[r] & M \\
S &
}
\]
constitute a $\Gamma$ equivariant family of closed $\Gamma$-invariant
proper subvarieties of $M$. Moreover the total space ${\mathcal Z}$ of
this family maps onto the $\Gamma$-invariant constructible subset 
$N_{i,k} \subset M$.

Let $x$ be a point in the open interior of $N_{i,k} \subset M$ and let
$Z \subset M$ be the Zariski closure of $\Gamma x$. Passing to a
subgroup of finite index in $\Gamma$ we can assume that $Z$ is
geometrically irreducible. Let $V \in S$ be such that $Z = Z_{V}$.
Then $x \in Z_{V}$ and by further localizing $S$ we may assume that
all the fibers of ${\mathcal Z} \to S$ are irreducible and of the same
dimension.

Let $W\subset {\mathcal Z}$ denote the set of all points which are contained
in two or more distinct fibers $Z_{V_1}$ and $Z_{V_2}$ (i.e. two fibers with
$Z_{V_1} \neq Z_{V_2}$). We
claim that this is a constructible subset. To see this,
look at the closed subvariety ${\mathcal Z}\times_{M} {\mathcal
Z} \subset {\mathcal Z} \times {\mathcal
Z}$ and the subset 
$I$ of $S\times S$ given by the condition
\[
I := \left\{ (V_{1},V_{2}) \in S\times S \left| 
\text{ 
such that $Z_{V_{1}} = Z_{V_{2}}$ as subvarieties in
$M$ 
}
\right.
\right\}.
\]
Note that $I \subset S\times S$ is a constructible subset (see this by
taking a compactification of $M$, compactifying the family ${\mathcal
Z} \to S$ and then using Chow schemes) and so the preimage ${\mathcal
I}$ of $I$ in ${\mathcal Z}\times {\mathcal Z}$ is also constructible.
The subset $W$ is the projection of ${\mathcal Z}\times_{M} {\mathcal
Z} - {\mathcal I}$ on one of the factors ${\mathcal Z}$, so $W$ is
constructible.

Next we claim that $W$ does not contain our original point $x$ thought
of as a point in the fiber $V\in S$. For if it did, this would mean
that there was a distinct $Z_{V'}\neq Z_{V}$ containing $x$, and then $x$
would be contained in the $\Gamma$-invariant set $Z_{V'}\cap Z_{V}$;
but this latter set is a proper subset of $Z_{V}$, contradicting the
fact (by definition of our family) that $\Gamma x$ is Zariski-dense in
$Z_{V}$.

It follows that there is an open set of the total space ${\mathcal Z}$
which does not meet $W$. Note that it is clear from the definition
that $W$ is the inverse image of a constructible subset in $M$; thus
this subset does not contain the generic point of $M$ so there is a
closed set $C$ such that $W$ is contained in the preimage of $C$ in
${\mathcal Z}$.  Since $W$ is a $\Gamma$-invariant set, we may replace
$C$ here by the intersection of all of its translates so we can assume
that $C$ is $\Gamma$-invariant as well. Finally then we can throw $C$
out of $M$ (i.e. replace $M$ by $M-C$ in the whole discussion) so we
may assume that $W$ is empty. Note that the new $M$ will no longer be
affine but only quasi-affine. This does not affect the rest of the
argument though since the existence of meromorphic functions can be
detected on opens.

By taking a compactification of $M$ and looking at the Chow scheme of
subvarieties of this compactification,
we can replace our family by a family indexed by a new base scheme $S$
where each
fiber (considered as a subset of $M$) occurs exactly once.

Now by the above reduction the morphism ${\mathcal Z}\rightarrow M$ is
injective on points; also its image hits a generic point of $M$. Thus
there is a largest open subset of $M$ over which this is an
isomorphism and we can replace $M$ by this open subset (which is
$\Gamma$-invariant). Hence we obtain a $\Gamma$-invariant fibration
$M\rightarrow S$. Now any meromorphic function on $S$ pulls back to
give a $\Gamma$-invariant meromorphic function on $M$. This
essentially proves the theorem. The only problem we need to address is
that in the construction of the family ${\mathcal Z} \to S$ we had to
pass to a finite index subgroup of $\Gamma$ and so the function just
constructed may be invariant only under a subgroup of finite index of
$\Gamma$.  This however is easily remedied - by taking the different
invariant polynomials in the Galois translates of our meromorphic
function we obtain a meromorphic function invariant by the full
$\Gamma$.  \eop

\

\medskip

\noindent
In view of Theorem~\ref{theo-alternative} we need to find effective
criteria for the non-existence of invariant rational functions on an
affine variety. One such criterion is given by the following lemma:

\begin{lem} \label{lem-functions} Let $\Gamma$ be a finitely presented
discrete group.
\begin{description}
\item[{\bf (alg)}] Suppose $M$ is a reduced irreducible scheme of
finite type over ${\mathbb C}$.  Suppose that $\Gamma$ acts on $M$ by
algebraic automorphisms and let $y\in 
M$ be a point in the smooth locus of $M$, fixed by the action, so that
$\Gamma$ acts linearly on the tangent space $T_{y}M$.  Let $G\subset
GL(T_{y}M)$ be the Zariski closure of $\op{im}[\Gamma\to
GL(T_{y}M)]$. Assume that $G$ acts on $T_{y}M$ with an open orbit,
and that the connected component $G^o$ of $G$ has no nontrivial
characters. Then there is no $\Gamma$-invariant rational
function on $M$. 
\item[{\bf (an)}]
Suppose that $N$ is an
irreducible analytic space on which 
$\Gamma$ acts by analytic automorphisms. 
Suppose that $y\in N$ is a point in the smooth locus of $N$, fixed by
the action, 
so that $\Gamma$ acts linearly on the tangent space $T_{y}N$.
Let $G\subset GL(T_{y}N)$ be the Zariski closure of $\op{im}[\Gamma\to
GL(T_{y}N)]$. Assume that $G$ acts on $T_{y}N$ with an open orbit,
and that the connected
component $G^o$ of $G$ has no nontrivial characters. 
Then there is no $\Gamma$-invariant analytic-meromorphic function on
$N$
\end{description}
\end{lem}
{\bf Proof.} Clearly the statement of the lemma is insensitive to
passing to a finite index subgroup of $\Gamma$ and so  we may
assume that $G = G^{o}$. We will only give the proof in the algebraic
case. The analytic case is completely analogous.

Suppose that $h$ is such a
function, and write the germ of $h$ at $y$ as $f/g$ with $f,g \in
{\mathcal O}_{M,y}$
relatively prime. Then for any $\gamma \in \Gamma$,
\[
\frac{f}{g} = h = \gamma ^{\ast}h = \frac{\gamma ^{\ast}f}{\gamma ^{\ast}g},
\]
which implies that there is a unit $u(\gamma )\in {\mathcal
O}_{M,y}^{\times}$ with
\[
\begin{split}
\gamma ^{\ast}f & = u(\gamma )f, \\
\gamma ^{\ast}g & = u(\gamma )g.
\end{split}
\]
Note that $\gamma \mapsto u(\gamma )$ is a cocycle for $\Gamma$
acting on
the multiplicative group of units ${\mathcal O}_{M,y}^{\times}$.  
In particular we get that
the value $u(\gamma )(y)$ is a character of $\Gamma$.

Without loss of generality we may assume that $h$ is not an
invertible function at $o$ (otherwise subtract the constant function
with the same value at $o$) hence we may assume that one of $f$ or $g$
has a nontrivial leading term of some degree. We may suppose that $f$
has such (otherwise replace $h$ by $h^{-1}$). Let $f_m$ be the leading
term of $f$.  Note that $f_m$ (of degree $m$)---modulo higher order
terms---is an element of $\op{Sym}^m T^{\vee}(M)_y$. The action of
$\Gamma$ on this leading term factors through the group $G$.  Our
previous formula gives
\[
\gamma ^{\ast}(f_m) = u(\gamma )(y) f_m.
\]
This shows that $\gamma \mapsto u(\gamma )(y)$ comes from a character of $G$;
but by assumption there are no such characters 
Therefore we get $\gamma ^{\ast}(f_m)=f_m$, so $f_m$ is a $G$-invariant
homogeneous
form of degree $m$.

In particular, we can think of $f_m$ as a $G$-invariant polynomial function on
the tangent space $T_{y}M$. This contradicts the supposed existence of an open
orbit in the action of $G$ on $T_{y}M$. The lemma is proven.
\eop

\

\medskip

\noindent
The essential consequence of Theorem~\ref{theo-alternative} that we
need is the following localization statement.

\begin{cor} \label{cor-localization}
Suppose that $y\in M$ is a point in the smooth locus of $M$, fixed by
the action, so that $\Gamma$ acts linearly on the tangent space
$T_{y}M$.  Let $G\subset GL(T_{y}M)$ be the Zariski closure of
$\op{im}[\Gamma\to GL(T_{y}M)]$. Assume that $G$ acts on $T_{y}M$
with an open orbit, and that the connected component $G^o$ of $G$ has
no nontrivial characters.  Then there exists a point $x\in M$ with
$\Gamma x$ Zariski-dense in $M$.
\end{cor}
{\bf Proof:} By Theorem~\ref{theo-alternative} we only have to rule
out the possibility that $M$ admits a nonconstant $\Gamma$-invariant
meromorphic function. This however is precisely the content of
Lemma~\ref{lem-functions}{\bf (alg)}. The Corollary is proven. \eop

\

Theorem~\ref{theo-alternative}, Lemma~\ref{lem-functions} and
Corollary~\ref{cor-localization} not only provide a convenient
localization criterion for the density of an action but also suggest
another geometric notion of `largeness' of the
$\Gamma$-orbits. Motivated by Theorem~\ref{theo-alternative}, we
define various degrees of analytic generic Zariski denseness ({\bf
AGZD} for short) as
follows:

\begin{defi} \label{defi-agzd} Suppose that a finitely generated group
$\Gamma$ acts by analytic automorphisms on an irreducible analytic
space $N$. Let $m : \Gamma\times N \to N$ be the action map.
\begin{itemize}
\item We say that the action $m$ is {\bf AGZD1} if there is no
$\Gamma$-invariant analytic meromorphic function $f$ on $N$.
\item We say that the action $m$ is {\bf AGZD2} if there is no pair 
$(U,f)$, where $U \subset N$ is a $\Gamma$-invariant analytically
Zariski dense open subset of $N$ and $f : U \to Z$ is a $\Gamma$
equivariant holomorphic map from $U$ to a complex analytic space $Z$
with $\dim Z < \dim N$.
\item We say that the action $m$ is {\bf AGZD3} if there is a point $x
\in N$ such that $m(\Gamma\times \{ x \}$ is analytically Zariski-dense in $N$.
\item We say that the action $m$ is {\bf AGZD4} if there is an 
analytically  Zariski dense open subset $U \subset N$ such that for every $x
\in U$ the $\Gamma$ orbit of $x$ is analytically Zariski-dense in $N$.
\end{itemize}
\end{defi}
\

\noindent
Clearly for an analytic action $m$ one has the implications:
\[
{\bf
AGZD4} \Rightarrow {\bf AGZD3} \Rightarrow {\bf AGZD2} \Rightarrow 
{\bf AGZD1},
\] 
but we don't think that the converse implications are
true. Similarly, it is clear that if $m$ is actually an algebraic
action, then {\bf AGZD1} implies that $m$ is Zariski-dense in the
algebraic sense of Theorem~\ref{theo-alternative}.

Suppose now that $B$ is a base scheme and that $p : M \to B$ is a
morphism equipped with a connection $\nabla$ (by which we mean a
stratification over the crystalline site of $S$
\cite{grothendieck-crystals}, \cite{simpson-moduli2}). For the
following definition it is not necessary to assume that $\nabla$ is
integrable. 

\begin{defi} \label{defi-gzd}
Suppose that $B$ and the generic geometric fiber of $M/B$
are irreducible. We say that $(p : M \to B,\nabla)$ is {\em
generically Zariski dense} (or {\bf GZD}) if there is no algebraic
meromorphic function $f$ on the total space $M$ which is invariant
under $\nabla$.
\end{defi}

If the connection $\nabla$ is integrable, then the corresponding analytic
family is associated to a local system of complex analytic spaces over
$B$, which in turn corresponds to the monodromy action $m$ of $\Gamma
:= \pi_{1}(B,o)$ on a fiber $M_{o}^{\op{an}}$. It is clear that if $m$
is {\bf AGZD1}, then $(p : M \to B,\nabla)$ is {\bf GZD}.

In particular, the {\bf AGZD1} property for $M_B(X_o,n)$ or equivalently
$M_{DR}(X_o,n)$ implies the algebraic generic Zariski-denseness
property {\bf GZD} for the Gauss-Manin connection $\nabla ^n_{DR}$.

\bigskip

Consider now a family of smooth projective connected curves $f : X \to B$ 
and let $o \in B$ be a base point. We will show that when the geometric
monodromy of $f : X \to B$ is of finite index in the mapping class
group or when $f$ comes from a Lefschetz pencil as in
Theorem~\ref{thm-lefschetz-pencil}, then the monodromy action of
$\pi_{1}(B,o)$ on $M_{B}(X_{o},n)$ is {\bf AGZD1}. In combination with
Theorem~\ref{theo-alternative} this fact yields statement (i) of
Theorems~\ref{thm-finite-index} and
\ref{thm-lefschetz-pencil}. As explained above this automatically
gives the analytic statement (ii) 
in both theorems. In fact, it
follows from the above considerations that $M_{DR}(X/B,n) \to B$
together with the non-abelian Gauss-Manin connection is {\bf GZD} in
the sense of Definition~\ref{defi-gzd}.

In view of all this it only remains to show that under the hypothesis
of Theorems~\ref{thm-finite-index} or \ref{thm-lefschetz-pencil} the
algebraic action 
\[
\monb{n} : \pi_{1}(B,o) \to
\op{Aut}(M_{B}(X_{o},n)).
\]
on the affine variety $M_{B}(X_{o},n)$ is  {\bf AGZD1}. (Note that
since $X_{o}$ is a smooth curve the variety 
$M_{B}(X_{o},n)$ is irreducible by \cite[Section~11]{simpson-moduli2}.)
In view of Corollary~\ref{cor-localization}, to achieve this we only
need to find a smooth point $\rho \in
M_{B}(X_{o},n)$ which is fixed by the monodromy
group $\monb{n}(\pi_{1}(B,o))$, and for
which the Zariski closure of $\monb{n}(\pi_{1}(B,o)) \subset
GL(T_{[\rho]}M_{B}(X_{o},n))$  acts on $T_{[rho]}M_{B}(X_{o},n)$ with
an open orbit and has a connected component of the identity which
admits no non-trivial characters.  

In the next section we describe a proposal for such a point $\rho$
which utilizes the Schr\"{o}dinger representation of a finite dihedral 
Heisenberg group. Later on, we will show in
Sections~\ref{sec-hyperelliptic} and \ref{sec-lefschetz} that the
open orbit property for the monodromy action on the tangent space at
$\rho$ holds, provided that the geometric monodromy of $f : X \to B$
is large enough.

We have stated the additional properties {\bf AGZD2-4} in order to
pose the question: which of these properties hold for families whose
monodromy has finite index in the mapping class group? For
(sufficiently ample) Lefschetz pencils?

\subsection{The Schr\"{o}dinger representation} \label{sec-Schrodinger}

Since $\rho$ is supposed to be fixed by the monodromy a  natural
choice would be to take $\rho$ to be the trivial representation of
$\pi_{1}(X_{o})$ in $\op{GL}(n,{\mathbb C})$. However the trivial
representation is a singular point of $M_{B}(X_{o},n)$ and so is
unsuitable for our 
purposes. On the other hand any representation
\[
\rho : \pi_{1}(X_{o}) \to \op{GL}(n,{\mathbb C})
\]
which has finite image will be fixed under some finite index
subgroup of $\monb{n}(\pi_{1}(B,o))$. Furthermore, the properties 
{\bf AGZD1-AGZD4} and  {\bf GZD} are obviously stable under passage to
a finite index subgroup of $\Gamma$. Hence we are free to replace $B$ by
any finite etale cover of $B$, and so it is enough to find a
finite representation $\rho$ which satisfies the open
orbit condition. 

In order to apply Corollary~\ref{cor-localization} we also need to
choose $\rho \in M_{B}(X_{o},n)$ 
to be a smooth point. This is  equivalent to choosing  $\rho$ to be an
irreducible representation.  

To construct such a representation we proceed as follows. Let $\bmu_{n}
\subset {\mathbb C}^{\times}$ be the group of all $n$-th roots of
unity.  Let $\widehat{{\mathbb Z}/n} := \op{Hom}_{{\mathbb
Z}}({\mathbb Z}/n,{\mathbb C}^{\times})$ denote the group of
characters of the cyclic group ${\mathbb Z}/n$. Consider the finite Heisenberg
group $H_{n}$. By definition $H_{n}$ is the central extension
\[
0 \to \bmu_{n} \to H_{n} \to {\mathbb Z}/n \times
\widehat{{\mathbb Z}/n} \to 0
\]
corresponding to the cocycle $e : ({\mathbb Z}/n \times
\widehat{{\mathbb Z}/n})^{2} \to \bmu_{n} \subset {\mathbb
C}^{\times}$, $e((a,\alpha),(a',\alpha')) = \alpha'(a)$. Explicitly
$H_{n}$ can be identified with the set $\bmu_{n}\times {\mathbb
Z}/n\times \widehat{{\mathbb Z}/n}$ with a group law given by
\begin{equation} \label{eq-h-grouplaw}
(\lambda; a, \alpha)\cdot (\lambda'; a', \alpha') =
(\lambda\lambda'\alpha'(a); a + a', \alpha\alpha').
\end{equation}
Let $\phi_{n} : H_{n} \to GL(V_{n})$ be the Schr\"{o}dinger
representation of $H_{n}$ \cite{mumford}. By definition $\phi_{n}$ is
the unique $n$-dimensional irreducible representation of $H_{n}$ which
has a tautological central character. One way to construct $\phi_{n}$
is to observe that the natural injective map
\[
\bmu_{n}\times \widehat{{\mathbb Z}/n} \hookrightarrow H_{n},\qquad
(\lambda, \alpha) 
\mapsto (\lambda; 0, \alpha)
\]
is a group monomorphism, i.e $H_{n}$ contains $\bmu_{n}\times \widehat{{\mathbb
Z}/n}$ as an abelian subgroup. Let ${\mathbb T}$ be the one dimensional
complex representation of $\bmu_{n}\times \widehat{{\mathbb
Z}/n}$ which corresponds to the pullback of the tautological character
of $\bmu_{n}$ under the projection $\bmu_{n}\times \widehat{{\mathbb
Z}/n} \to \bmu_{n}$. In other words ${\mathbb T} = ({\mathbb C},\tau)$
where $\tau : \bmu_{n}\times {\mathbb
Z}/n \to {\mathbb C}^{\times}$ is given by $\tau(\lambda,\alpha) =
\lambda$. In terms of ${\mathbb T}$ then we have 
\[
(V_{n},\phi_{n}) =
\op{Ind}_{\bmu_{n}\times \widehat{{\mathbb Z}/n}}^{H_{n}}({\mathbb
T}). 
\]
Explicitly we can identify $V_{n}$ with the vector space of all
complex valued functions on the finite set ${\mathbb Z}/n$ and the
action $\phi_{n}$ by the formula
\[
[\phi_{n}(\lambda; a, \alpha)f](x) = \lambda\alpha(x)f(x+a),
\]
for all $x \in {\mathbb Z}/n$, $f \in V_{n}$ and $(\lambda; a, \alpha)
\in H_{n}$. 

The irreducibility of the representation $\phi_{n}$ follows from
Frobenius reciprocity or directly by noticing that $V_{n}$ has a basis
consisting of the characteristic functions of the elements in
${\mathbb Z}/n$ and that the subgroup ${\mathbb Z}/n \subset H_{n}$
acts transitively on the elements of this basis.  In particular if we
compose $\phi_{n}$ with some surjective homomorphism 
$\pi_{1}(X _{o}) \twoheadrightarrow H_{n}$ we will get a representation
of $\pi_{1}(X_{o})$ in $GL(V_{n})$ which is irreducible and has finite
image. Unfortunately, it turns out (see
Remark~\ref{rem-duality}) that this representation can not be used
directly to obtain an open orbit action on the tangent space to 
$M_{B}(X_{o},n)$. However a slight modification of this representaion
does the job. The modification involves an extension of $H_{n}$ of
dihedral type which we proceed to describe.

Let $\bmu_{2} = \{ -1, +1 \} \subset {\mathbb C}^{\times}$ be the
group of square roots of one. The group $\bmu_{2}$ acts naturally on
${\mathbb Z}/n\times \widehat{{\mathbb Z}/n}$ as the inversion on both
factors. This action clearly preserves the cocycle defining the
Heisenberg central extension $0 \to \bmu_{n} \to H_{n} \to {\mathbb
Z}/n\times \widehat{{\mathbb Z}/n} \to 0$ and so we get a natural
action of $\bmu_{2}$ on $H_{n}$. Explicitly, if we think of $H_{n}$ as
the set $\bmu_{n}\times {\mathbb Z}/n \times \widehat{{\mathbb Z}/n}$
equipped with the group law \eqref{eq-h-grouplaw}, then an element
$\varepsilon \in \bmu_{2}$ acts on $H_{n}$ via $(\lambda; a, \alpha)
\mapsto (\lambda, \varepsilon a,\alpha^{\varepsilon})$. We define the
{\em dihedral Heisenberg group} as the semidirect product
\[
{\mathfrak D}H_{n} := \bmu_{2} \ltimes H_{n}
\]
for the above action. Thus ${\mathfrak D}H_{n}$ can be identified with the set 
$\bmu_{2}\times \bmu_{n}\times {\mathbb Z}/n \times \widehat{{\mathbb
Z}/n}$ with a group law given by
\begin{equation} \label{eq-dh-grouplaw}
(\varepsilon,\lambda,a,\alpha)\cdot (\varepsilon',\lambda',a',\alpha')
=
(\varepsilon\varepsilon',\lambda\lambda'\alpha^{'\varepsilon}(a), 
a+\varepsilon 
a', \alpha\alpha^{'\varepsilon}).
\end{equation}
In particular, for each $(\varepsilon,\lambda,a,\alpha) \in {\mathfrak
D}H_{n}$ we have
\[
\begin{split}
(\varepsilon,\lambda,a,\alpha) & = (1,\lambda,a,\alpha)\cdot
(\varepsilon,1,0,\bone) \\
& = (\varepsilon,1,0,\bone)\cdot (1,\lambda,\varepsilon
a,\alpha^{\varepsilon}),
\end{split}
\]
where $\bone : {\mathbb Z}/n \to {\mathbb C}^{\times}$ stands for the
trivial character, i.e. $\bone(a) = 1$ for all $a$.

Observe next that the Schr\"{o}dinger representation $\phi_{n} : H_{n}
\to GL(V_{n})$ extends naturally to a {\em dihedral Schr\"{o}dinger
representation} 
\[
{\mathfrak d}\phi_{n} : {\mathfrak D}H_{n} \to GL(V_{n})
\]
defined by
\[
({\mathfrak d}\phi_{n}(\varepsilon,\lambda,a,\alpha)f)(x) =
\lambda\alpha(x) f(\varepsilon(x+a))
\]
for all $f \in V_{n}$ and all $x \in {\mathbb Z}/n$. 

Recall that if $C$ is any smooth curve of genus $g$, then there
is a surjective homomorphism $\pi_{1}(C) \twoheadrightarrow
\boldsymbol{F}_{g}$ onto a free group of $g$ generators. This
homomorphism is obtained by moding $\pi_{1}(C)$ out by the normal
subgroup generated by the $a$-cycles for a standard basis in the first
homology of $C$. Note furthermore that ${\mathfrak D}H_{n}$ is generated
by the three
elements $(-1,1,0,\bone)$, $(1,1,1,\bone)$ and $(1,1,0,\alpha)$, where
$\alpha \in \widehat{{\mathbb Z}/n}$ is any generator. Hence we can find a
surjective homomorphism $\pi_{1}(C) \to \boldsymbol{F}_{g} \to
{\mathfrak D}H_{n}$ 
as long as $g \geq 3$.

By hypothesis the genus of $X_{o}$ is big enough and so we can find a
surjective homomorphism $\psi_{n} : \pi_{1}(X _{o}) \twoheadrightarrow
{\mathfrak D}H_{n}$. Let $\rho : \pi_{1}(X_{o}) \to \op{GL}(V_{n})$
denote the composition $\rho := \psi_{n}\circ {\mathfrak
d}\phi_{n}$. By construction $\rho$ is irreducible and so represents a
smooth point of the moduli of representations. Moreover by a standard
deformation theory argument (see e.g. \cite{repsfg}) we can identify
the Zariski tangent space $T_{[\rho]}M_{B}(X_{o},n)$ with the group
cohomology $H^{1}(\pi_{1}(X_{o}), \op{ad}(\rho))$, where
\[
\op{ad}(\rho) : \pi_{1}(X_{o}) \to \op{GL}(\op{End}(V_{n}))
\]
is the natural representation induced from $\rho$. Explicitly
$\op{ad}(\rho) = ({\mathfrak d}\phi_{n}^{\vee}\otimes {\mathfrak
d}\phi_{n})\circ \psi_{n}$ and 
since $\op{ad}(\phi_{n}) = \phi_{n}^{\vee}\otimes \phi_{n}$ has a
trivial central character we see that $\op{ad}(\rho)$ factors trough
the quotient group ${\mathfrak D}H_{n} \twoheadrightarrow
\bmu_{2}\ltimes ({\mathbb Z}/n\times \widehat{{\mathbb Z}/n})$.  

It is not hard to calculate $H^{1}(\pi_{1}(X_{o}), \op{ad}(\rho))$ in
terms of geometric data on the curve $X_{o}$. The action of $\bmu_{2}$
on ${\mathbb Z}/n\times \widehat{{\mathbb Z}/n}$ induces an obvious
action (inversion on both factors) of $\bmu_{2}$ on the group characters 
\[
\op{Hom}({\mathbb
Z}/n\times \widehat{{\mathbb Z}/n},{\mathbb C}^{\times}) =
\widehat{{\mathbb Z}/n}\times {\mathbb Z}/n.
\] 
Now for each orbit $u
\in (\widehat{{\mathbb Z}/n}\times {\mathbb Z}/n)/\bmu_{2}$ we get an
irreducible representation $W_{u}$ of \linebreak $\bmu_{2}\ltimes ({\mathbb
Z}/n\times \widehat{{\mathbb Z}/n})$. The representation $W_{(\bone,0)}$
corresponding to the trivial character is just the trivial one
dimensional representation of $\bmu_{2}\ltimes ({\mathbb
Z}/n\times \widehat{{\mathbb Z}/n})$. For any other orbit $u$ we have
that $u = \{ \chi, \chi^{-1} \}$ for some non-trivial character $\chi
\in \widehat{{\mathbb Z}/n}\times {\mathbb Z}/n$ and so $W_{u}$ is the
representation of $\bmu_{2}\ltimes ({\mathbb
Z}/n\times \widehat{{\mathbb Z}/n})$ induced from the one dimensional
representation $({\mathbb C},\chi)$ of ${\mathbb
Z}/n\times \widehat{{\mathbb Z}/n}$. Thus $W_{u} = ({\mathbb
C},\chi)\oplus ({\mathbb C},\chi^{-1})$ as a representation of ${\mathbb
Z}/n\times \widehat{{\mathbb Z}/n}$ and the generator of $\bmu_{2}$
acts by switching the two summands. In particular $W_{u}$ is a two
dimensional irreducible (we are assuming that $n$ is odd here)  
representation of $\bmu_{2}\ltimes ({\mathbb
Z}/n\times \widehat{{\mathbb Z}/n})$. 

With this notation we have

\begin{lem} \label{lem-tangent-space} The tangent space to
$M_{B}(X_{o},n)$ at the dihedral Schr\"{o}dinger representation $\rho$
is given by  
\begin{equation} \label{eq-tangent-space}
T_{[\rho]}M_{B}(X_{o},n) = H^{1}(X_{o},\op{ad}(\rho)) =
\bigoplus_{u
\in (\widehat{{\mathbb Z}/n}\times {\mathbb Z}/n)/\bmu_{2}} H^{1}(X_{o},{\mathbb
W}_{u}) 
\end{equation}
where ${\mathbb W}_{u}$ is the local system on $X_{o}$
corresponding to the representation
\[
\pi_{1}(X_{o}) \stackrel{\psi_{n}}{\to} {\mathfrak D}H_{n} \to
\bmu_{2}\ltimes ({\mathbb Z}/n\times \widehat{{\mathbb Z}/n}) \to
\op{GL}(W_{u}).
\]
\end{lem}
{\bf Proof.} Note that a representation $\kappa : {\mathbb Z}/n\times
\widehat{{\mathbb Z}/n} \to GL(V)$ of the abelian group ${\mathbb
Z}/n\times \widehat{{\mathbb Z}/n}$ will extend to representation of
the `dihedral' group $\bmu_{2}\ltimes ({\mathbb Z}/n\times
\widehat{{\mathbb Z}/n})$ if and only if $\kappa$ is
self-dual. Furthermore, each self-dual representation $\kappa$ has a
canonical dihedral extension:
\[
{\mathfrak d}\kappa : \bmu_{2}\ltimes ({\mathbb Z}/n\times
\widehat{{\mathbb Z}/n}) \to \op{GL}(V),
\]
in which $\bmu_{2}$ acts as the self-daulity automorphism of
$V$. Concretely if we decompose $(V,\kappa)$ into a direct sum of
characters of ${\mathbb Z}/n\times \widehat{{\mathbb Z}/n}$, then the
self-duality of $V$ will identify the multiplicity space of each
character $\chi$ with the multiplicity space of the character
$\chi^{-1}$. In particular the multiplicity spaces in the character
decomposition of $(V,\kappa)$ depend not on the individual characters
but rather on the $\bmu_{2}$-orbits $u \in (\widehat{{\mathbb
Z}/n}\times {\mathbb Z}/n)/\bmu_{2}$. Hence $(V,{\mathfrak d}\kappa)$
decomposes as
\[
(V,{\mathfrak d}\kappa) = \bigoplus_{u \in (\widehat{{\mathbb
Z}/n}\times {\mathbb Z}/n)/\bmu_{2}} W_{u}\otimes M_{u},
\]
where $M_{u}$ denotes the multiplicity space of a character $\chi \in
u$ in $(V,\kappa)$.

Consider now the representation $\op{ad}(\phi_{n})$ of $H_{n}$. Since
it has a trivial central character, it factors through a representation
\[
\op{ad}(\phi_{n})^{\op{ab}} : {\mathbb Z}/n\times \widehat{{\mathbb
Z}/n} \to \op{GL}(\op{End}(V_{n})).
\]
This is the {\em abelian part} of
$\op{ad}(\phi_{n})$.  The representation $\op{ad}(\phi_{n})^{\op{ab}}$
is self-dual by construction and so admits a canonical dihedral
extension 
\[
{\mathfrak d}\op{ad}(\phi_{n})^{\op{ab}} : \bmu_{2}\ltimes ({\mathbb
Z}/n\times \widehat{{\mathbb Z}/n}) \to \op{GL}(\op{End}(V_{n})).
\]
This dihedral extension 
fits in the commutative diagram
\[
\xymatrix{
\pi_{1}(X_{o}) \ar[d]_-{\psi_{n}} \ar[r]^-{\op{ad}(\rho)} & 
\op{GL}(\op{End}(V_{n})) \\
{\mathfrak D}H_{n} \ar[r] &  \bmu_{2}\ltimes ({\mathbb
Z}/n\times \widehat{{\mathbb Z}/n})
\ar[u]_-{{\mathfrak d}\op{ad}(\phi_{n})^{\op{ab}}} 
}
\]
and so understanding $\op{ad}(\rho)$ is equivalent to understanding
$\op{ad}(\phi_{n})^{\op{ab}}$. But $\op{ad}(\phi_{n})^{\op{ab}}$ is
just the regular representation of the abelian group ${\mathbb
Z}/n\times \widehat{{\mathbb Z}/n}$. To see this note first that  $\dim
\op{End}(V_{n}) = n^{2} = \dim {\mathbb C}[{\mathbb
Z}/n\times \widehat{{\mathbb Z}/n}]$. Now since every irreducible
representation of ${\mathbb
Z}/n\times \widehat{{\mathbb Z}/n}$ occurs in  ${\mathbb C}[{\mathbb
Z}/n\times \widehat{{\mathbb Z}/n}]$ with multiplicity one we need only
to check that for every character $\chi  : {\mathbb Z}/n\times
\widehat{{\mathbb Z}/n} \to {\mathbb C}^{\times}$ we have 
\[
\dim \op{Hom}_{({\mathbb Z}/n\times
\widehat{{\mathbb Z}/n})-\op{mod}}(\chi,\op{ad}(\phi_{n})^{\op{ab}})
\geq 1.
\]
But the group of characters of ${\mathbb Z}/n\times
\widehat{{\mathbb Z}/n}$ is naturally isomorphic to $\widehat{{\mathbb
Z}/n}\times {\mathbb Z}/n$ and so each character $\chi$ as above is
given by a pair $(\xi,x) \in \widehat{{\mathbb
Z}/n}\times {\mathbb Z}/n$ via the formula $\chi(a,\alpha) =
\xi(a)\cdot \alpha(x)$. Therefore we only need to show that for any
pair $(\xi,x)$ there exists a non zero element $A_{(\xi,x)} \in
\op{End}(V_{n}) = V_{n}^{\vee}\otimes V_{n}$ so that
\[
(a,\alpha)A_{(\xi,x)} = \xi(a)\alpha(x)A_{(\xi,x)} 
\]
for all $(a,\alpha) \in {\mathbb Z}/n\times \widehat{{\mathbb Z}/n}$.

To construct the element $A_{(\xi,x)}$ recall that the vector space
$V_{n} = {\mathbb C}[{\mathbb Z}/n]$ has a natural basis $\{ e_{0},
e_{1}, \ldots, e_{n-1} \}$ consisting of characteristic functions of
elements of ${\mathbb Z}/n$, i.e. $e_{i}(j) := \delta_{ij}$. Let $\{
e_{0}^{\vee}, e_{1}^{\vee}, \ldots, e_{n-1}^{\vee} \}$ be the dual
basis of $V_{n}^{\vee}$. Then  in terms of the 
basis $\{ e_{i}^{\vee}\otimes e_{j} \}_{i,j =
0}^{n-1}$ of $V_{n}^{\vee}\otimes V_{n}$ the representation
$\op{ad}(\phi_{n})^{\op{ab}}$ is given by the formula
\[
[\op{ad}(\phi_{n})^{\op{ab}}(a,\alpha)](e_{i}^{\vee}\otimes e_{j}) =
\alpha(j-i)e_{i-a}^{\vee}\otimes e_{j-a}.
\]
In view of this we may take
\[
A_{(\xi,x)} := \sum_{i=0}^{n-1} \xi(i) e_{i}^{\vee}\otimes e_{i+x}
\]
which is obviously a non-zero eigenvector corresponding to the
character $\chi = (\xi,x)$.

This shows that 
\[
(\op{End}(V_{n}),\op{ad}(\phi_{n})^{\op{ab}}) = {\mathbb C}[{\mathbb
Z}/n\times \widehat{{\mathbb Z}/n}] = \oplus_{\chi \in
\widehat{{\mathbb Z}/n}\times {\mathbb Z}/n} ({\mathbb C}, \chi),
\]
and so the lemma is proven. \eop

\

\bigskip

Let us now go back to the problem of checking whether
$\op{mon}(\pi_{1}(B,o))$ has a dense orbit on $M_{B}(X_{o},n)$. As
mentioned above, the fact that $\rho$ has a finite image implies that
the conjugacy class of $\rho$ will be fixed by some finite-index
subgroup of $\op{mon}(\pi_{1}(B,o))$. In particular, applying
Corollary~\ref{cor-localization} to this subgroup, it follows that in
order to show the existence of a dense $\op{mon}(\pi_{1}(B,o))$-orbit
on $M_{B}(X_{o},n)$, it suffices to check the following two items:
\begin{description}
\item[(i)] The Zariski closure $G$ of 
\[
\op{im}[\op{mon}(\pi_{1}(B,o)) \to
GL(H^{1}(X_{o},\op{ad}(\rho)))]
\] 
in $GL(H^{1}(X_{o},\op{ad}(\rho)))$
has an open orbit on $H^{1}(X_{o},\op{ad}(\rho))$. 
\item[(ii)] The identity component $G^{o}$ of $G$ does not have any
non-trivial characters.
\end{description}

\

\medskip

\noindent
Condition {\bf (ii)} follows easily from the isomorphism
(\ref{eq-tangent-space}):

\begin{lem} \label{lem-condition-ii}
The identity component $G^{o}$ of the Zariski closure of
\[\op{im}[\op{mon}(\pi_{1}(B,o)) \to 
GL(H^{1}(X_{o},\op{ad}(\rho)))]
\] 
in $GL(H^{1}(X_{o},\op{ad}(\rho)))$
has no non-trivial characters.
\end{lem}
{\bf Proof.} Indeed, let $p: Y_{o} \to X_{o}$ be the dihedral
Galois cover of $X_{o}$ corresponding to the surjection 
\begin{equation} \label{eq-abelian-quotient}
\pi_{1}(X_{o}) \stackrel{\psi_{n}}{\to} {\mathfrak D}H_{n} \to
\bmu_{2}\ltimes ({\mathbb Z}/n\times 
\widehat{{\mathbb Z}/n}).
\end{equation}
Then $Y_{o}$ is a smooth connected curve and the pushforward of the
trivial one dimensional local system ${\mathbb C}_{Y_{o}}$ via $p$ is 
precisely
\[
p_{*}{\mathbb C}_{Y_{o}} = \bigoplus_{u \in (\widehat{{\mathbb
Z}/n}\times {\mathbb Z}/n)/\bmu_{2}} {\mathbb W}_{u}.
\]
Thus we can identify $H^{1}(X_{o},\op{ad}(\rho))$ with
$H^{1}(Y_{o},{\mathbb C})$. 

Next observe that without a loss of generality we may assume that the
family $f : X \to B$ has an algebraic section $\sigma : B \to
X$. Indeed, since $X$ is quasi-projective a generic hyperplane section
on $X$ will be a multisection of $f$. But replacing $B$ by an
\'{e}tale cover of a Zariski open set of $B$ will replace
$\op{mon}(\pi_{1}(B,o))$ by a subgroup of finite index. In particular
such a replacement will not affect the property that 
 $\op{mon}(\pi_{1}(B,o))$ fixes the
conjugacy class of $\rho$. In fact, by taking another \'{e}tale cover
if necessary we can ensure that not only the conjugacy class of $\rho$
is fixed under $\op{mon}(\pi_{1}(B,o))$ but that the actual
representation
\[
\rho : \pi_{1}(X_{o},\sigma(o)) \to GL(V_{n})
\]
remains fixed under $\op{mon}(\pi_{1}(B,o))$. Indeed, to achieve this
we only need to pass to the finite index subgroup of the monodromy
which preserves the kernel of $\rho$ and acts trivially on the finite
group $\op{im}(\rho) = {\mathfrak D}H_{n}$.

Assume that we are in this situation. Then $\rho$ lifts to a well
defined representation of $\pi_{1}(B,o)\ltimes_{\op{mon}}
\pi_{1}(X_{o},\sigma(o))$ and so defines a ${\mathfrak D}H_{n}$-cover
of $X$. In particular under the identification
$H^{1}(X_{o},\op{ad}(\rho)) \cong H^{1}(Y_{o},{\mathbb C})$ the action
of $\op{mon}(\pi_{1}(B,o))$ on $H^{1}(X_{o},\op{ad}(\rho))$ becomes
just the monodromy action of $\pi_{1}(B,o)$ on $H^{1}(Y_{o},{\mathbb
C})$ corresponding to the family of curves $Y \to B$. But by Deligne's
semisimplicity theorem \cite[Corollaire~4.2.9]{deligne-hodge2} the
monodromy action on the middle dimensional cohomology of any smooth
projective family over a quasi-projective base has a semisimple
Zariski closure. Thus $G^{o}$ must be a connected semisimple algebraic
group and so has no non-trivial characters. The lemma is proven. \eop

\

\bigskip

\noindent
In order to check that condition {\bf (i)} is satisfied we need to
make sure that the monodromy group of the family $f : X \to B$ is
sufficently large.

Before we explain how this is achieved  we need to introduce some notation.
On the way we will also rephrase the condition {\bf (i)} in a slightly
more general context.

Let $f : X \to B$ be a smooth family of connected curves of genus $g
\geq 3$. Let as before $o \in B$ be a fixed base point and let 
\[
\op{mon} : \pi_{1}(B,o) \to \op{Map}(X_{o}) \subset \op{Out}(\pi_{1}(X_{o}))
\]
be the corresponding geometric monodromy representation. 

If $f : X \to B$ has a
holomorphic section $\sigma : B \to X$  the representation 
\[
\op{mon} : \pi_{1}(B,o) \to \op{Map}(X_{o})
\] 
can be lifted to a
geometric monodromy representation respecting the base point:
\[
\op{mon}^{\sigma} : \pi_{1}(B,o) \to \op{Map}^{1}(X_{o}) \subset 
\op{Aut}(\pi_{1}(X_{o},\sigma(o))).
\]
Here $\op{Map}^{1}(X_{o})$ denotes the mapping class group of the once
punctured surface $X_{o} - \{ \sigma(o) \}$.

Fix a finite abelian group $A$ and let ${\mathfrak D}A :=
\bmu_{2}\ltimes A$ denote the standard dihedral extension of
$A$ in which the generator $(-1) \in \bmu_{2}$ acts as the inversion on
$A$. Fix a surjective homomorphism $\pi_{1}(X_{o},\sigma(o))
\twoheadrightarrow {\mathfrak D}A$ and let $p^{{\mathfrak D}A} :
Y^{{\mathfrak D}A}_{o} \to X_{o}$ be the 
corresponding Galois cover. Let $\op{Map}^{1}(X_{o},{\mathfrak D}A)$ be
the group of 
$p^{{\mathfrak D}A}$-liftable mapping classes, i.e. 
\[
\op{Map}^{1}(X_{o},{\mathfrak D}A) = \left\{ \varphi \in
\op{Map}^{1}(X_{o}) \left| 
\begin{minipage}[c]{3.2in}
$\varphi$ preserves $\ker[\pi_{1}(X_{o},\sigma(o))
\twoheadrightarrow {\mathfrak D}A]$  and $\varphi$ induces the
identity on ${\mathfrak D}A$ 
\end{minipage}
\right. \right\}
\]

Clearly $\op{Map}^{1}(X_{o},{\mathfrak D}A) \subset
\op{Map}^{1}(X_{o})$ is of finite index and consists precisely of the
mapping classes on $X_{o}$ which lift to mapping classes on
$Y_{o}^{{\mathfrak D}A}$. Furthermore, since by definition each
$\varphi \in \op{Map}^{1}(X_{o},{\mathfrak D}A)$ induces the identity
on ${\mathfrak D}A$ it follows that any lift $\tilde{\varphi} \in
\op{Map}(Y_{o}^{{\mathfrak D}A})$ of $\varphi$ commutes with the
action of ${\mathfrak D}A$ on $Y_{o}^{{\mathfrak D}A}$, and that
${\mathfrak D}A$ acts transitively on the set of all such lifts. Thus,
if we define $\op{LMap}^{1}(X_{o},{\mathfrak D}A)\subset
\op{Map}^{1}(Y_{o}^{{\mathfrak D}A})$ to be the subgroup consisting of
all lifts of elements in $\op{Map}^{1}(X_{o})$ we see that
$\op{LMap}^{1}(X_{o},{\mathfrak D}A)$ fits in a short exact sequence
of groups
\[
1 \to {\mathfrak D}A \to \op{LMap}^{1}(X_{o},{\mathfrak D}A) \to
\op{Map}^{1}(X_{o},{\mathfrak D}A) \to 1. 
\]
\

\noindent
Assume now that $\op{mon}^{\sigma}(\pi_{1}(B,o)) \subset
\op{Map}^{1}(X_{o},{\mathfrak D}A)$. In particular
$\op{mon}^{\sigma}(\pi_{1}(B,o))$ preserves $\pi_{1}(Y_{o}^{{\mathfrak
D}A})$ and 
so we get a short exact sequence of groups 
\[
1 \to \pi_{1}(B,o)\times_{\op{mon}^{\sigma}} \pi_{1}(Y_{o}^{{\mathfrak
D}A}) \to 
\pi_{1}(X_{o},\sigma(o)) \to {\mathfrak D}A \to 1. 
\]
Let $Y^{{\mathfrak D}A} \to X$ denote the ${\mathfrak D}A$-Galois
cover of $X$ corresponding to the homomorphism $\pi_{1}(X,\sigma(o))
\twoheadrightarrow {\mathfrak D}A$. By construction
$\pi_{1}(Y^{{\mathfrak D}A}) \cong
\pi_{1}(B,o)\times_{\op{mon}^{\sigma}} \pi_{1}(Y_{o}^{{\mathfrak
D}A})$ and the corresponding monodromy representation
\[
\op{mon}_{{\mathfrak D}A} : \pi_{1}(B,o) \to
\op{Map}(Y_{o}^{{\mathfrak D}A}) \subset 
\op{Out}(\pi_{1}(Y_{o}^{{\mathfrak D}A})) 
\]
lands in $\op{LMap}^{1}(X_{o},{\mathfrak D}A)$.

Furthermore, note that from the viewpoint of the density properties we
are interseted in, the conditions that $f : X \to B$ has a section and
that $\op{mon}(\pi_{1}(B,o)) \subset \op{Map}^{1}(X_{o},{\mathfrak
D}A)$ are harmless. Indeed, as explained in the proof of
Lemma~\ref{lem-condition-ii}, if $f : X \to B$ is an arbitrary smooth
projective family of curves with $B$ smooth and connected and $X$
quasi-projective, then we can always find a Zariski open set $U
\subset B$ containing the point $o \in B$, and a finite \'{e}tale
cover $(B',o') \to (U,o)$, so that the pulled-back family $X\times_{B}
B' \to B'$ has a holomorphic section, and a geometric monodromy which
is contained in $\op{Map}^{1}(X_{o},{\mathfrak D}A)$. Since
$\pi_{1}(U,o) \twoheadrightarrow \pi_{1}(B,o)$ is surjective and
$\pi_{1}(B',o') \subset \pi_{1}(U,o)$ is a subgroup of finite index,
it follows that the geometric monodromy $\op{mon}(\pi_{1}(B',o'))$ of
the family $X\times_{B} B' \to B'$ is a subgroup of finite index in
$\op{mon}(\pi_{1}(B,o))$. In particular, any density statement we can make for
the action of $\pi_{1}(B,o)$ will be equivalent to the corresponding
density statement for the action of $\pi_{1}(B',o')$.

The previous reasoning also shows that for any smooth family of curves
\linebreak $f : X \to B$, such that $B$ is smooth and $X$ is quasi-projective,
and any surjective homomorphism $\pi_{1}(X_{o},\sigma(o))
\twoheadrightarrow {\mathfrak D}A$, there is an appropriate $(B',o')
\to (B,o)$ and an ${\mathfrak D}A$-Galois cover $Y^{{\mathfrak D}A}
\to X\times_{B} B'$, so that:
\begin{itemize}
\item The image of the 
monodromy representation $\op{mon}_{{\mathfrak D}A} : \pi_{1}(B',o') \to
\op{Map}(Y_{o}^{{\mathfrak D}A})$ is contained in
$\op{LMap}^{1}(X_{o},{\mathfrak D}A)$; 
\item The natural map $\op{mon}_{{\mathfrak D}A}(\pi_{1}(B',o')) \to
\op{mon}(\pi_{1}(B,o))$ has finite kernel and cokernel.
\end{itemize}

Motivated by the discussion in Sections~\ref{sec-open-orbits} and
\ref{sec-Schrodinger} we make the following definition:

\begin{defi} \label{defi-good} Let $A$ be a finite abelian group. A
pair $(f : X \to B, \pi_{1}(X_{o}) \twoheadrightarrow {\mathfrak D}A)$ is
called {\em good} if the Zariski closure of 
\[
\op{im}\left[\pi_{1}(B',o')  
\xrightarrow{\op{mon}_{{\mathfrak D}A}}
\op{LMap}^{1}(X_{o},{\mathfrak D}A) \subset \op{Map}(Y_{o}^{{\mathfrak D}A}) \to
\op{Sp}(H_{1}(Y_{o}^{{\mathfrak D}A}, {\mathbb Z}))\right]
\]
in $\op{Sp}(H_{1}(Y_{o}^{{\mathfrak D}A}, {\mathbb C}))$ acts on
$H_{1}(Y_{o}^{{\mathfrak D}A}, {\mathbb C})$ with an open orbit.
\end{defi}

Clearly now, the condition {\bf (i)} is equivalent to the statement
that if $A = {\mathbb
Z}/n\times \widehat{{\mathbb Z}/n}$ and if 
the homomorphism $\pi_{1}(X_{o}) 
 \to {\mathfrak D}A$ is induced from a surjective
homomorphism 
\[
\pi_{1}(X_{o}) \to {\mathfrak D}H_{n},
\] 
then the pair
$(f : X \to B, \pi_{1}(X_{o}) \twoheadrightarrow {\mathfrak D}A)$ is good.

In other words we need to find geometric restrictions on a family $f :
X \to B$ and a \linebreak homomorphism $\pi_{1}(X_{o}) \twoheadrightarrow
{\mathfrak D}A$, which will guarantee that the pair \linebreak $(f : X \to B,
\pi_{1}(X_{o}) \twoheadrightarrow {\mathfrak D}A) $ is good.

As a first approximation one has to understand the image of
$\op{LMap}^{1}(X_{o},{\mathfrak D}A)$ into the symplectic group
$\op{Sp}(H_{1}(Y_{o}^{{\mathfrak D}A},{\mathbb Z})$. In the next
section we will analyze the hyperelliptic part of this image for a
suitably chosen surjection $\pi_{1}(X_{o}) \twoheadrightarrow
{\mathfrak D}A$.

\

\bigskip

\begin{rem} \label{rem-duality} Satisfying condition {\bf (i)} is a
somewhat subtle task. In a preliminary version of this paper we
attempted to work with a representation $\rho : \pi_{1}(X_{o}) \to
\op{GL}(V_{n})$ which comes from a choice of a surjective homomorphism
$\pi_{1}(X_{o}) \twoheadrightarrow H_{n}$ onto the Heisenberg group
$H_{n}$ rather than its dihedral extension ${\mathfrak D}H_{n}$. This
representation $\rho$ is also irreducible and so gives a smooth point
in $M_{B}(X_{o},n)$ which satisfies condition {\bf (ii)}. Furthermore
the image of $\op{LMap}(X_{o},A)$ into the corresponding symplectic
group $\op{Sp}(H_{1}(Y_{o}^{A},{\mathbb Z})$ was described explicitly
by Looijenga \cite[Theorem~2.5]{looijenga}. Unfortunately the
self-duality pairing on the representation
$\op{ad}(\phi_{n})^{\op{ab}}$ gives rise to a quadratic function on
$H^{1}(X_{o},\op{ad}(\rho))$ which will be preserved by all elements of
the geometric monodromy and so we can not hope that $G$ will have an
open orbit for this choice of $\rho$. Replacing $\rho$ by a
representation coming from a surjection onto the dihedral Heisenberg
group repairs this problem as we will see below. However this changes
the setup and forces us to work with the
two-dimensional dihedral representations $W_{u}$ instead of the
characters of $A$. In particular this setup lies beyond the scope of
Looijenga's analysis in \cite{looijenga} and forces us to look for a
description of the image of $\op{LMap}^{1}(X_{o},{\mathfrak D}A)$ into
the symplectic group $\op{Sp}(H_{1}(Y_{o}^{{\mathfrak D}A},{\mathbb
Z}))$ based on first principles only. 

Remarkably enough, it turns out that such a concrete description is
possible and that it leads to a stronger result which uses only
hyperelliptic mapping classes to obtain an open orbit. However note
that that we need to assume that $n$ is odd in the explicit argument.
\end{rem}

\section{Proofs of the main theorems} \label{sec-proofs}
\setcounter{figure}{0}

In this section we prove
Theorems~\ref{thm-finite-index} and \ref{thm-lefschetz-pencil}.

\subsection{The case of a hyperelliptic monodromy} \label{sec-hyperelliptic}

In this section $f : X \to B$ will denote a smooth family of all
hyperelliptic curves of genus $g$. For us the hyperelliptic curves
will be represented as branched double covers of ${\mathbb P}^{1}$
having $2g + 2$ branch points. From that point of view it is natural
to take $B$ to be the configuration space $\op{{\mathbf Conf}}_{2g+2}$
of $2g+2$ distinct points in ${\mathbb P}^{1}$. However, it is well
known \cite{mess} that no universal hyperelliptic family exists on
that space. This can be remedied by either passing to an unramified
double cover of $\op{{\mathbf Conf}}_{2g+2}$ (see \cite[Remark~4]{mess}) or,
alternatively, by taking an open subfamily of
$\op{{\mathbf Conf}}_{2g+2}$. We take the second approach since it is
better suited for 
our purposes.

Concretely, we take  $B$ to be the configuration space
of $2g+2$ distinct points in the affine line
${\mathbb C} = {\mathbb P}^{1}-\{\infty
\}$. To see that the universal hyperelliptic family on $B$ exists, we
only need to show that the incidence divisor 
\[
\Sigma := \left\{\left.
(\{b_{1},\ldots,b_{2g+2}\},x) \in 
B\times {\mathbb P}^{1}\right| x \in \{b_{1},\ldots,b_{2g+2}\}\right\} \subset
B\times {\mathbb P}^{1}
\]
corresponds to a section $\sigma$ in the line bundle $p_{{\mathbb
P}^{1}}^{*}{\mathcal O}(2g +2)$ on $B\times {\mathbb P}^{1}$. Indeed,
if this is the case we have a divisor $\sigma(B)$ in the total space
$\op{tot}(p_{{\mathbb P}^{1}}^{*}{\mathcal O}(2g +2))$ of the line
bundle $p_{{\mathbb P}^{1}}^{*}{\mathcal O}(2g +2)$ and so $X$ can be
constructed simply as the preimage of $\sigma(B)$ in
$\op{tot}(p_{{\mathbb P}^{1}}^{*}{\mathcal O}(g+1))$ under the natural
squaring map.

To see that $\Sigma$ is in the linear system $|p_{{\mathbb
P}^{1}}^{*}{\mathcal O}(2g +2)|$ one can argue as follows. Fix an
affine coordinate $z$ on ${\mathbb C} = {\mathbb
P}^{1}-\{\infty\}$. For any integer $k > 0$ let $S_{k} :=
H^{0}({\mathbb P}^{1},{\mathcal O}_{{\mathbb P}^{1}}(k))$. On the
product $S_{k} \times {\mathbb P}^{1}$ we have the line bundle
$p_{{\mathbb P}^{1}}^{*}{\mathcal O}(k)$. Moreover the direct image
$p_{S_{k}*}p_{{\mathbb P}^{1}}^{*}{\mathcal O}(k)$ is a vector bundle
of rank $k+1$ on $S_{k}$ which is canonically isomorphic to ${\mathcal
O}_{S_{k}}\otimes S_{k}$ and so has a tautological section
corresponding to $\op{id}_{S_{k}}$. Let $\sigma_{k}$ denote the
corresponding section of $p_{{\mathbb
P}^{1}}^{*}{\mathcal O}(k)$. By construction, the divisor of
$\sigma_{k}$ consists of all pairs $(s,x) \in S_{k}\times {\mathbb
P}^{1}$ with $s(x) = 0$. Consider now the subspace 
$E_{k} \subset S_{k}$ of all monic
polynomials in $z$ of degree $k$. Since by definition $B$ can be
identified with the
open subset of the afiine 
subspace $E_{2g+2} \subset S_{2g+2}$ consisting of monic
polynomials with simple zeros, we can take $\sigma =
\sigma_{2g+2|B}$. 

Therefore we have constructed a universal double cover
\[
\xymatrix{
X \ar[rr]^-{\nu}\ar[dr]_-{f} & & B\times {\mathbb P}^{1}
\ar[dl]^-{p_{B}} \\
& B &
}
\]
over the configuration space $B$. Note that the fundamental group of
$B$ is the braid group $B_{2g+2}$ on $2g+2$ strands and that the
monodromy homomorphism for $f : X \to B$ can be interpreted as the
standard surjection from $B_{2g+2}$ onto the hyperelliptic mapping
class group. 

Fix as base point $o \in B$ the double cover 
$\nu_{o} : X_{o} \to {\mathbb P}^{1}$ with branch points $1,2,\ldots ,
2g+2 \in {\mathbb R} \subset {\mathbb P}^{1}$ on
the real axis. Make branch cuts $C_{i}$ on ${\mathbb P}^{1}$ from
$2i-1$ to $2i$ along the real axis. Topologically, the surface $X_{o}$
is obtained (see Figure~\ref{fig-glue} below) 
by gluing together two copies of the sliced-up ${\mathbb
P}^{1}$ along the rims of the branch cuts.

\begin{figure}[!ht]
\begin{center}
\psfrag{P1}[c][c][0.9][0]{{${\mathbb P}^{1}$}}
\psfrag{C1}[c][c][0.7][0]{{$C_{1}$}}
\psfrag{C2}[c][c][0.7][0]{{$C_{2}$}}
\psfrag{C(g+1)}[c][c][0.7][0]{{$C_{g+1}$}}
\psfrag{cut}[c][c][0.7][0]{{cut}}
\psfrag{glue}[c][c][0.7][0]{{glue}}
\epsfig{file=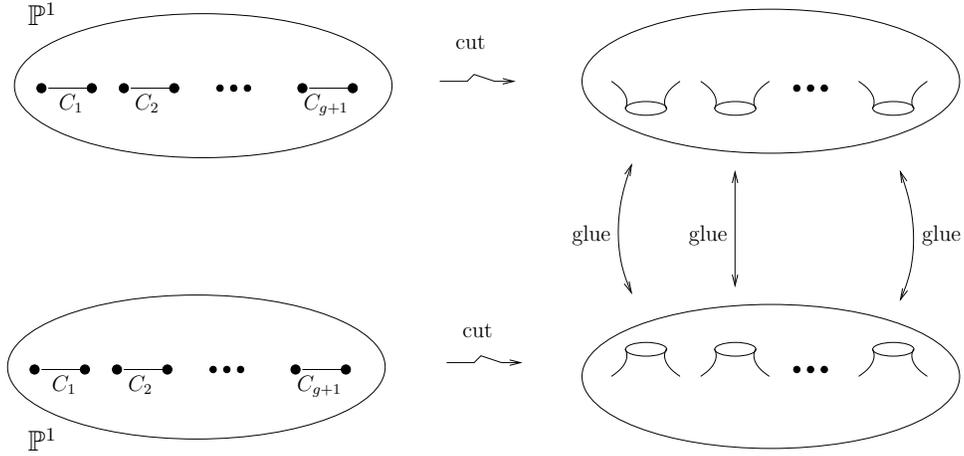,width=5in} 
\end{center}
\caption{Gluing two ${\mathbb P}^{1}$ sheets into a hyperelliptic curve.}
\label{fig-glue} 
\end{figure}

\

For concretness we label the two ${\mathbb P}^{1}$-sheets as the {\em
upper} and the {\em lower } sheet of $X_{o}$. The covering map
$\nu_{o} : X_{o} \to {\mathbb P}^{1}$ projects each sheet onto
${\mathbb P}^{1}$ and the corresponding covering involution
$\imath_{o} : X_{o} \to X_{o}$ interchganges the two sheets.

While working with loops on $X_{o}$ (either as representatives of
elements in $\pi_{1}(X_{o})$ or as circles determining Dehn twists on
$X_{o}$) it will be convenient to describe these loops in terms of
their $\nu_{o}$-images in ${\mathbb P}^{1}$. We will only look at
loops which do not pass through the branch points and which are
transversal to the boundary circles of our sheets. Every such loop $L$
projects via $\nu_{o}$ onto a simple closed path in ${\mathbb P}^{1}$
which does not pass through any branch point and is transversal to the
branch cuts. Therefore specifying a loop $L$ in $X_{o}$ is the same
thing as specifying a simple closed path in ${\mathbb P}^{1}$ (not
passing through the branch points and transversal to the branch cuts),
together with a labeling at each point indicating whether it is on the
upper or lower sheet and such that this labeling changes upon crossing
a branch cut. To avoid introducing additional notation we will write
$L$ both for the loop in $X_{o}$ and for the corresponding labeled
path in ${\mathbb P}^{1}$. This will not create any confusion since it
will always be clear from the context which incarnation of $L$ we have
in mind.

\begin{figure}[!ht]
\begin{center}
\psfrag{P1}{{${\mathbb P}^{1}$}}
\psfrag{Xo}{{$X_{o}$}}
\psfrag{nuo}{{$\nu_{o}$}}
\psfrag{a1}{{$a_{1}$}}
\psfrag{a2}{{$a_{2}$}}
\psfrag{a3}{{$a_{3}$}}
\psfrag{a4}{{$a_{4}$}}
\psfrag{a5}{{$a_{5}$}}
\psfrag{a6}{{$a_{6}$}}
\psfrag{a2g+1}{{$a_{2g+1}$}}
\psfrag{a2g+2}{{$a_{2g+2}$}}
\psfrag{cut}{{cut}}
\psfrag{glue}{{glue}}
\epsfig{file=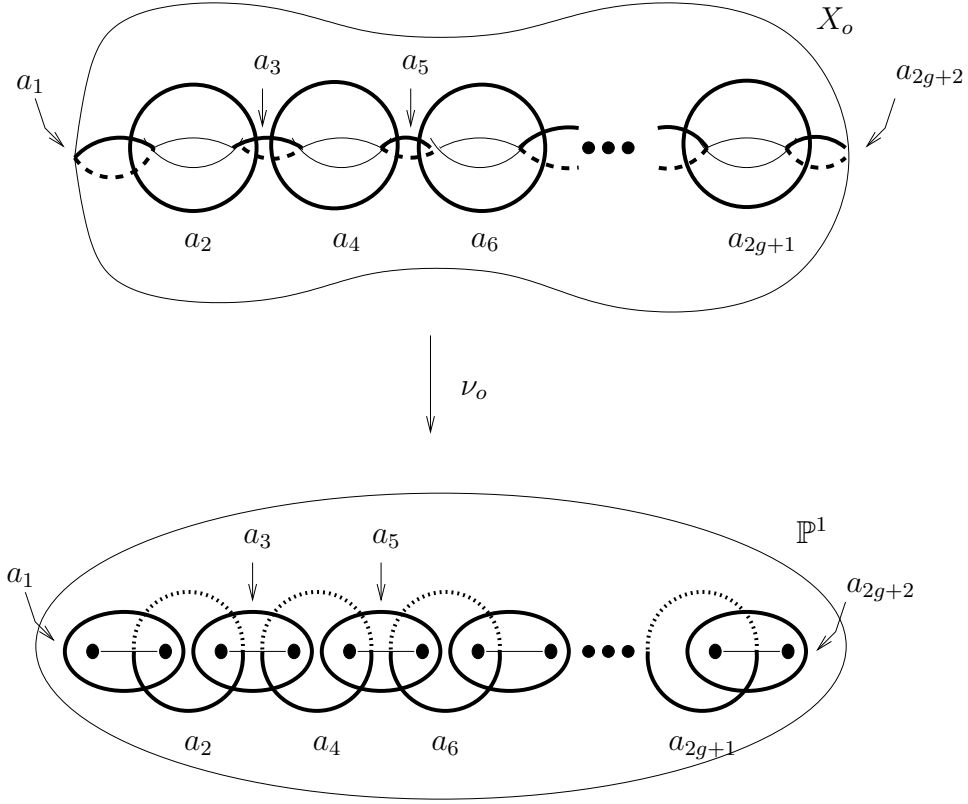,width=5in} 
\end{center}
\caption{Dehn twists generating $\Delta(X_{o})$.}
\label{fig-heloops} 
\end{figure}

The image of the geometric monodromy representation
\[
\op{mon} : \pi_{1}(B,o) \to \op{Map}(X_{o})
\]
for the family $f : X \to B$, is the full hyperelliptic mapping class
group 
\[
\Delta(X_{o}) = \{ \phi \in \op{Map}(X_{o}) |
\phi\imath_{o}\phi^{-1} = \imath_{o} \}.
\] 
It is generated by the
right handed Dehn twists along the sequence of 
loops $a_{1}, \ldots, a_{2g+1}$ depicted on Figure~\ref{fig-heloops}. 
In this picture we use the convention that for  paths
in ${\mathbb P}^{1}$ the solid pieces are on the upper sheet and the
dotted pieces are on the lower sheet.

\

\noindent
Note that the above Dehn twists define a surjective homomorphism from
the braid group $B_{2g+2}$ on $2g + 2$ strands to the hyperelliptic
mapping class group $\Delta(X_{o})$. Indeed, by definition $B_{2g+2}$
can be presented as
\[
B_{2g+2} = \left\langle t_{1}, t_{2}, \ldots, t_{2g+1} \left| \;
\begin{minipage}[c]{2.1in}
\[
\begin{split}
t_{i}t_{j} & = t_{j}t_{i}, \text{ for $|i-j| \geq 2$} \\
t_{i}t_{i+1}t_{i} & = t_{i+1}t_{i}t_{i+1}
\end{split}
\]
\end{minipage}
\right. \right\rangle
\]
and so the assignment $t_{i} \mapsto \text{(Dehn twist along
$a_{i}$)}$ induces a (necessarily surjective) group homomorphism
$\kappa_{g} : B_{2g+2} \to \Delta(X_{o})$.

\

Fix a positive odd integer $n$. To fix notation, choose  a primitive
$n$-th root of unity  $\gamma \in \bmu_{n}$ and let $\alpha : {\mathbb
Z}/n \to {\mathbb C}^{\times}$, $k \mapsto \gamma^{k}$ be the corresponding
character. For future use we denote the corresponding standard
generators of ${\mathfrak D}H_{n}$ as follows
\[
\bsigma := (-1,1,0,\bone), \qquad
\ba := (1,1,1,\bone), \qquad
\balpha := (1,1,0,\alpha).
\]
We also write $\be := (1,1,0,\bone)$ for the identity element in
${\mathfrak D}H_{n}$. 

With this notation we are now ready to define the base representation
$\rho \in M_{B}(X_{o},n)$ at which we will be checking the open orbit
condition {\bf (i)} from section~\ref{sec-Schrodinger}. For this we
only need to exhibit a surjective homomorphism $\psi_{n} :
\pi_{1}(X_{o}) \to {\mathfrak D}H_{n}$. We define $\psi_{n}$ to be
trivial on the complement of the branch cuts on both sheets and we
postulate that the passing transformations $P_{i} \in {\mathfrak
D}H_{n}$ corresponding to going through the branch cut $C_{i}$ should be:
\[
P_{1} = \bsigma, \quad
P_{2} = \ba \bsigma, \quad
P_{3} = \balpha \bsigma, \quad
P_{4} = \ba, \quad
P_{5} = \balpha, \quad
P_{i} = \be, \quad \text{for all $i \geq 6$}.
\]
Assume that $g\geq 6$, so there are at
least two branch cuts with passing transformation equal to the
identity. 

Consider next any two element orbit $u = \{ \chi, \chi^{-1} \} \in
(\widehat{{\mathbb Z}/n}\times {\mathbb Z}/n)/\bmu_{2}$ and let $W_u$
be the corresponding $2$-dimensional irreducible representation of
${\mathfrak D}H_{n}$. As we saw in section~\ref{sec-Schrodinger}, the
action of ${\mathfrak D}H_{n}$ on $W_{u}$ factors through
$\bmu_{2}\ltimes ({\mathbb Z}/n\times \widehat{{\mathbb Z}/n})$. Choose
a basis $\{ v_{+}, v_{-} \}$ of $W_{u}$ consisting of eigenvectors
for the ${\mathbb
Z}/n\times \widehat{{\mathbb Z}/n}$ action. To fix notation assume
that $v_{+}$ corresponds to the character $\chi$ and that $v_{-}$
corresponds to the character $\chi^{-1}$. If we write the character
$\chi$ as $\chi = (\balpha^{b},\ba^{c})$ for some 
integers $b$ and $c$, then in the basis $\{ v_{+}, v_{-} \}$ 
the representation $W_{u}$ is given by associating 
\[
\bsigma \mapsto R, \quad
\ba \mapsto P^{b}, \quad
\balpha  \mapsto P^{c},
\] 
where $P$ and $R$ are the $2\times 2$ matrices:
\[
P := \begin{pmatrix} \gamma & 0 \\ 0 & \gamma^{-1} \end{pmatrix}, \qquad
R := \begin{pmatrix} 0 & 1 \\ 1 & 0 \end{pmatrix}.
\]
This gives the matrices for the action of passing transformations on
the local system ${\mathbb W}_{u}$:
\begin{equation} \label{eq-passing}
P_{1} = R, \qquad
P_{2} = P^b R, \qquad
P_{3} = P^c R, \qquad
P_{4} = P^b, \qquad
P_{5} = P^c.
\end{equation}
and the rest are equal to the identity matrix $I \in
\op{GL}_{2}({\mathbb C})$.

Note that by our assumption on $n$ it follows that $P$ is of odd order
so these matrices are never equal to $-I$.

\begin{figure}[!ht]
\begin{center}
\psfrag{P1}{{${\mathbb P}^{1}$}}
\psfrag{Xo}{{$X_{o}$}}
\psfrag{nuo}{{$\nu_{o}$}}
\psfrag{L12}{{$L_{12}$}}
\psfrag{L23}{{$L_{23}$}}
\psfrag{L47}{{$L_{47}$}}
\psfrag{L68}{{$L_{68}$}}
\epsfig{file=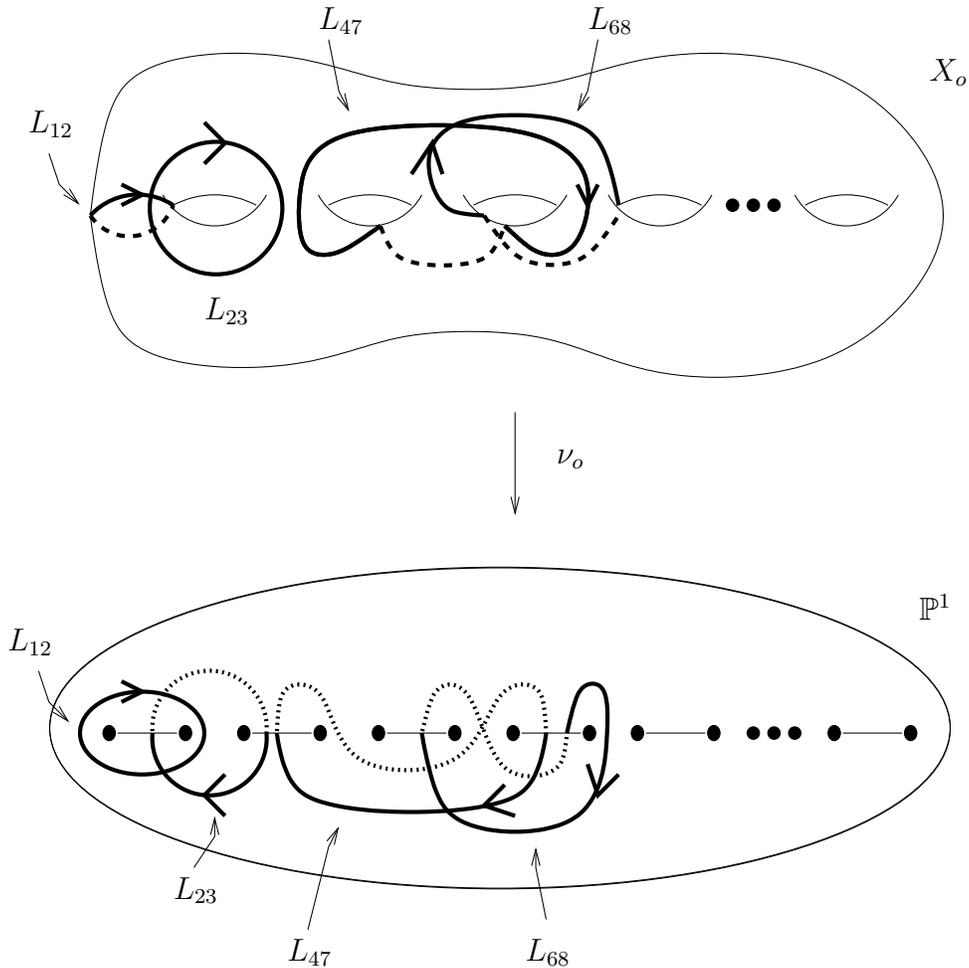,width=5in} 
\end{center}
\caption{Some $L_{ij}$'s.}
\label{fig-generators} 
\end{figure}

For any $1 \leq i < j \leq 2g+2$ let $L_{ij}$ denote the loop which goes
around the branch points $i$ and $j$, passing under any other branch
points which are in between on the real axis, and going in the
clockwise direction. (Some sample loops $L_{ij}$ are illustrated on
Figure~\ref{fig-generators}.) Assume that the lower part of the curve
is on the upper sheet, and let $M_{ij} \in \op{GL}_{2}({\mathbb C})
\cong \op{GL}(W_{u})$ denote the monodromy transformation around
$L_{ij}$.

Observe next that the representation $W_{u}$ is self-dual: in our
basis $\{v_{+},v_{-}\}$ the invariant pairing $Q : W_{u}\otimes W_{u}
\to {\mathbb C}$ is given by
\[
Q\left(\begin{pmatrix} r \\ s \end{pmatrix},  \begin{pmatrix} r' \\ s'
\end{pmatrix}\right) := rs' + r's.
\]
The pairing $Q$ induces an intersection pairing on
$H_{1}(X_{o},{\mathbb W}_{u})$ and an isomorphism of homology and cohomology
$H_{1}(X_{o},{\mathbb W}_{u}) \cong H^{1}(X_{o},{\mathbb W}_{u})$,
both compatible with the monodromy action of the hyperelliptic mapping
class group $\pi_{1}(B,o) \cong \Delta(X_{o})$. Thus it suffices to
calculate the monodromy action on the homology  $H_{1}(X_{o},{\mathbb
W}_{u})$. 

We will represent the elements in $H_{1}(X_{o},{\mathbb W}_{u})$ by
loop-like chains. Similarly to ordinary loops, a looplike chain with
coefficients in $W_{u}$ is given by:
\begin{itemize}
\item an oriented simple closed path in ${\mathbb P}^{1}$ (not passing
through the branch points and transversal to the branch cut), together
with a labeling at each point indicating whether it is on the upper or
the lower sheet, and such that the labeling changes upon crossing of a
branch cut;
\item a specification at each point of the path of a vector in
$W_{u}$, such that upon crossing a branch cut $C_{i}$ from the upper
to the lower sheet this vector is
modified by the corresponding passing matrix $P_{i}$.
\end{itemize}
A looplike chain will typically be denoted by $vL$, where $L$ is the
loop and $v$ is the corresponding vector in $W_{u}$.

Our convention for the intersection pairing on $H_{1}(X_{o},{\mathbb
W}_{u})$  will be that we intersect the looplike chains by
intersecting the underlying loops according to the right-hand rule (a
first path intersects positively a second path which, to him is coming
from the right) and at each intersection point we pair the
corresponding elements in $W_{u}$ via $Q$.

The homology group $H_{1}(X_{o},{\mathbb W}_{u})$ contains elements of the
form $v_{ij}L_{ij}$, where $v_{ij} \in W_{u}$ is any vector which is
invariant under $M_{ij}$.  We can distinguish three cases: 
\begin{itemize}
\item[$\lozenge$] if $M_{ij}$ is the
identity then $v_{ij}$ can be any vector so there is a two-dimensional
space of such cycles; 
\item[$\lozenge$] if $M_{ij}$ is a reflection (of the form $P^kR$)
then $v_{ij}$ is of the form $(v + M_{ij}v)$ for any sufficiently
general vector in
$W_{u}$, and  in fact we may take
$v = v_{+}$;
\item[$\lozenge$] if $M_{ij}$ is a rotation of the form $P^k$ for $k$
different from $0$ modulo the order of $P$, then there are no nonzero
cycles of this form.
\end{itemize}
Using these elements we can now prove the following:

\begin{lem} \label{lem-span}
The cycles of the form $v_{ij}L_{ij}$ span the homology
$H_{1}(X_{o},{\mathbb W}_{u})$ of the hyperelliptic curve with
coefficients in ${\mathbb W}_u$.
\end{lem}
{\bf Proof.}  We need to separate into cases depending on $b$ and
$c$. Assume first of all that $b$ and $c$ are different and different
from $0$ (modulo $n$).  For the purposes of this lemma, we can apply a
braid transformation to arrange things so that the passing matrices
are (in order starting with $P_{1}$):
\[
P^b, \; 
P^c, \; 
P^bR, \; 
P^cR, \; 
R, \; 
I,\ldots .
\]
Now we consider an element of the homology of the hyperelliptic curve
with coefficients in ${\mathbb W}_u$.  It can be moved to a cycle
supported over 
the real axis, necessarily on the interval between $1$ and $2g+2$.  Look
first at the interval $[2g+1, 2g+2]$.  By subtracting off a cycle of the
form $v_{2g+1,2g+2}L_{2g+1,2g+2}$ we obtain a cycle which is zero on the
upper sheet along the interval in question (note that the vector
$v_{2g+1,2g+2}$ can be arbitrary).  Now the cycle condition at the point
$2g+2$ implies that the cycle is also zero on the lower sheet. We get to
a cycle supported on the interval $[1,2g+1]$. Continuing this way by
induction we get to a cycle supported on the interval $[1,10]$ (the last
nontrivial passing matrix is $P_{5}$).

The monodromy transformation $M_{9,10}$ is the identity and so again
we can take $v_{9,10}$ to be an arbitrary vector. By subtracting a
multiple of $v_{9,10}L_{9,10}$ we get to a 
cycle supported in $[1,9]$.

Next look at the monodromy matrices
\[
M_{2,9} = RP^b\;\; \text{and} \;\; M_{4,9} = RP^c.
\]
Thus we get that up to scalars
\[
v_{2,9}= (1+RP^b)v_+ = 
\begin{pmatrix}
1   \\
\gamma^b  
\end{pmatrix} 
\]
and 
\[
v_{4,9}= (1+RP^c)v_+ = 
\begin{pmatrix}
1   \\
\gamma ^c  
\end{pmatrix}.
\]
The determinant of the matrix with these two vectors in the columns is
$\gamma^c -\gamma^b$ which is nonzero under our assumption that $b$
is different from $c$.  Therefore by subtracting off an appropriate
combination of the cycles $v_{2,9}L_{2,9}$ and $v_{4,9}L_{4,9}$ we
obtain a cycle which is zero in the interval $[8,9]$, in other words it
is supported on $[1,8]$.  Again the monodromy matrix $M_{7,8}$ is the
identity so we can subtract off a vector of the form $v_{7,8}$ to get
a cycle supported in $[1,7]$.

We repeat the argument above using 
\[
\begin{split}
v_{2,7} & = (1 + P^bRP^b)v_+ = v_+, \; \text{and} \\
v_{4,7} & = (1+P^bRP^c)v_+= v_+ + \gamma ^{c-b}v_-
\end{split}
\]
which are linearly independent. This, combined with subtracting off a
$v_{5,6}L_{5,6}$, gets us to a cycle supported in $[1,5]$; and repeating
again the same argument with $v_{2,5}$ and $v_{4,5}$ we get to a cycle
supported in $[1,3]$.  On the other hand the monodromy transformation
$M_{1,2}$ is trivial so by subtracting off a cycle of the form
$v_{1,2}L_{1,2}$ with $v_{1,2}$ arbitrary, we get to a cycle supported
in $[2,3]$, which must be a multiple of $v_{2,3}L_{2,3}$ so we are done.
One can note in passing that this last cycle must automatically be
zero since the monodromy $M_{2,3}$ doesn't have any fixed vectors but
this is not really important.
This completes the proof in the case $b\neq c$.  

Suppose now we are in the
case $b=c$ (modulo $n$). Then $b$ and $c$ are nonzero.  Thus we can
repeat the same argument as above but arranging things so that the
passing matrices are (in order)
\[
P^b, \; 
P^c, \; 
I,\;
P^bR, \; 
P^cR, \; 
R, \; 
I,\ldots .
\]
In this case the same argument as before (but using vectors such as
$v_{4,11}$ and $v_{6,11}$ etc.)  allows us to get to a cycle
supported on $[1, 5]$.  Now the first monodromy matrices are the
identities:
\[
M_{1,2} = M_{2,3}= M_{3,4} = I
\]
so we can subtract off cycles of the form $v_{1,2}L_{1,2}$,
$v_{2,3}L_{2,3}$, and $v_{3,4}L_{3,4}$ to get to a cycle supported on
$[4,5]$ and again we are done.  This completes the proof of the lemma.
\eop

\

\bigskip

\noindent
For the remainder of the argument we return to the labeling
\eqref{eq-passing} for the order of the branch cuts.

\begin{figure}[!ht]
\begin{center}
\psfrag{P1}[c][c][0.9][0]{{${\mathbb P}^{1}$}}
\psfrag{L34}[c][c][0.9][0]{{$L_{34}$}}
\psfrag{D23}[c][c][0.9][0]{{$D_{23}$}}
\epsfig{file=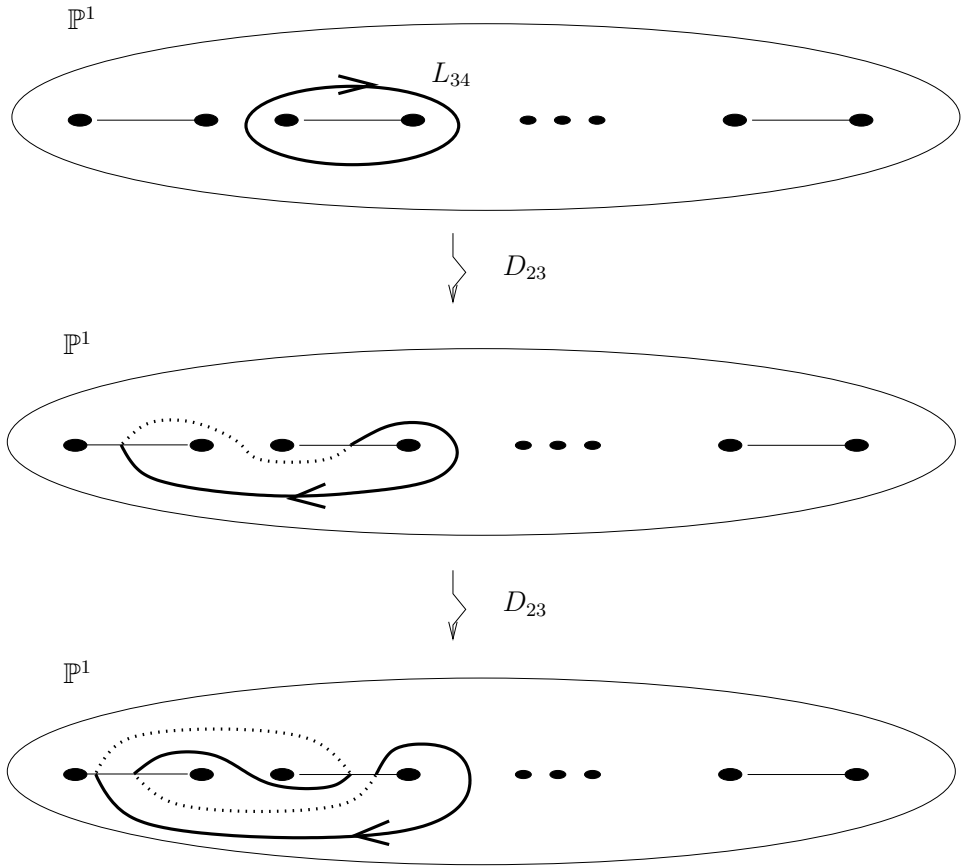,width=5in} 
\end{center}
\caption{The action of $D_{23}^{2}$ on the loop $L_{34}$.}
\label{fig-dehn} 
\end{figure}

Now let $\boldsymbol{t}_{ij} \in \op{Map}(X_{o})$ 
denote the right handed Dehn twist along the loop $L_{ij}$. For the
three types of behavior of $M_{ij}$ we have:
\begin{itemize}
\item[$\lozenge$] if $M_{ij}$ is the
identity, then $\boldsymbol{t}_{ij} \in \Delta(X_{o})$  and acts on the local
system $R^{1}f_{*}{\mathbb W}_{u}$. In particular
$\boldsymbol{t}_{ij}$ maps to a 
well defined element $D_{ij} \in \op{Sp}(H_{1}(X_{o},{\mathbb W}_{u}))$; 
\item[$\lozenge$] if $M_{ij}$ is a reflection, then
$\boldsymbol{t}_{ij}^2 \in \Delta(X_{o})$ and acts on the local system
$R^{1}f_{*}{\mathbb W}_{u}$. In particular $\boldsymbol{t}_{ij}^{2}$
maps to a well defined element $D_{ij}^{2} \in \op{Sp}(H_{1}(X_{o},{\mathbb
W}_{u}))$ (Figure~\ref{fig-dehn} illustrates the typical action of $D_{ij}^2$); 
\item[$\lozenge$] if $M_{ij}$ is a rotation then we don't consider the Dehn
twist.  
\end{itemize}
\

\

\noindent
For uniformity we will always consider $D_{ij}^2$. We now have
the following lemma:

\newpage

\begin{lem} \label{lem-irreducible} The subgroup of
$\op{Sp}(H_{1}(X_{o},{\mathbb W}_{u}))$ generated by the elements
\[
\left\{ D_{ij}^2 \left| \begin{minipage}[c]{1.3in} $1 \leq i <
j \leq 2g+2$ \\ $M_{ij} \neq$ rotation 
\end{minipage}
\right.
\right\}
\]
acts irreducibly on the complex vector space
$H_{1}(X_{o},{\mathbb W}_{u})$.
\end{lem}
{\bf Proof.} 
Consider the elements of the the form
$D_{ij}^2 -1$ in the group algebra of $\op{Sp}(H_{1}(X_{o},{\mathbb
W}_{u}))$. By general principles this algebra 
is semisimple, so
it suffices to find a vector $w_0$ such that the subspace generated by
the action of the $D_{ij}^2-1$ starting with $w_0$, spans the whole
$H_{1}(X_{o},{\mathbb
W}_{u})$.  We will start with
\[
w_0:= v_{2,11}L_{2,11}
\]
and show that using the $D_{ij}^2-1$ we can get to any vector
$v_{ij}L_{ij}$.  In view of Lemma~\ref{lem-span} this will complete the
proof of the irreducibility.

Note that $P_{6}$ is the identity and $P_{1}$ is a reflection, so
$M_{2,11}=R$ is a reflection. Thus $v_{2,11}$ is of the form $v_+
+Rv_+ = v_+ + v_-$.

Using $D_{12}^2-1$ we get to $(v_{+}+Rv_{+})L_{12}$ (note that we allow
ourselves to multiply by a factor for example $\frac{1}{2}$ or
$-\frac{1}{2}$ when we 
say this).

Now one of $b$ or $c$ is different from zero modulo the order of
$P$. Assume for example that $b$ is different from zero.  Then the
monodromy transformation $M_{2,7}$ is $P^{-b}R$ so using $D_{27}^2-1$
we get to
\[
(1+P^{-b}R)(1+R)v_{+}L_{27}.
\]
Applying again $D_{12}^2-1$ we get back to
\[
(1+P^{-b}R)(1+R)v_{+}L_{12}.
\]
In the other case where $b$ is zero but $c$ nonzero we could use
$D_{29}^{2} - 1$
and get to
\[
(1+P^{-c}R)(1+R)v_{+}L_{12}.
\]
In the first case, note that the image of $(1+R)$ is not contained in
the kernel of $(1+ P^{-b}R)$, and the image of $(1+P^{-b}R)$ is
linearly independent from the image of $(1+R)$, so with the vector we
obtained previously we obtain both vectors $v_+L_{12}$ and
$v_-L_{12}$. The same holds in the second case where we used the
matrix $P^c$.

A similar argument gets us to any of the vectors $vL_{34}$ and $vL_{56}$.  

Next, using again the fact that one of the $P^bR$ or $P^cR$ is
different from $R$, and using the appropriate transformation
$D_{4,11}^2-1$ or $D_{6,11}^2-1$ as well as $D_{2,11}^2-1$ and
following by $D_{11,12}^2-1$, we get to any vector of the form
$vL_{11,12}$.  Now using the Dehn twists for $i,j$ with $11 \leq i < j
\leq 2g+2$ we obtain all of the vectors of the form $vL_{ij}$ for $11
\leq i < j \leq 2g$.

From these using the Dehn twists $D_{i,j}^2-1$ for $i< 11$ and $j\geq
11$ we get to all vectors of the form $v_{ij}L_{ij}$ for $i< 11$ and
$j\geq 11$.

Similarly we get to all vectors of the form $v_{ij}L_{ij}$ when the
monodromy transformations $M_{ij}$ are reflections, by using the Dehn
twist $D_{ij}^2-1$ on a vector $v_{kl}L_{kl}$ where one of $k$ or $l$
is either $i$ or $j$ and where $[k,l]$ is a branch cut on which the
passing matrix is a reflection (note that such $k,l$ always exist when
the monodromy $M_{ij}$ is a reflection).  Here from above we have
already gotten to the $v_{kl}L_{kl}$ with $v_{kl}$ arbitrary.

This argument also works to obtain $v_{ij}L_{ij}$ whenever $1\leq i <
j \leq 6$ is an exceptional case where the monodromy is the identity
due to a special equality of the form $b=c$ or $b=0$ or $c=0$.

The only cycles which remain to be obtained are the $v_{ij}L_{ij}$ for
$7\leq i < j \leq 10$.  We first obtain $v_{78}L_{78}$.  To do this,
note that two among the three matrices $R$, $P^bR$ and $P^cR$ are
different, and for appropriate choices of $i$ and $j$ corresponding to
these two, chosen among $2$, $4$ and $6$, we have that the images of
the rank one matrices $(1+M_{i,7})$ and $(1+M_{j,7})$ are linearly
independent and span our two dimensional space (this is similar to the
argument used in Lemma~\ref{lem-span}).  Thus applying the Dehn twist
$D_{78}^2-1$ to the vectors $v_{i,7}$ and $v_{j,7}$ we span a two
dimensional space so we can get to any vector of the form
$v_{78}L_{78}$ with $v_{78}$ arbitrary.  The same argument yields any
vector of the form $v_{9,10}L_{9,10}$ with $v_{9,10}$
arbitrary. Finally, in the exceptional case where $M_{89}$ is the
identity (this is when $b=c$), using its Dehn twist we get to the
vectors of the form $v_{89}L_{89}$.

This completes the proof of the lemma. \eop

\

\bigskip

Let now ${\mathfrak g}$ denote the Lie algebra of the monodromy group
acting on the representation ${\mathbb C}^{4g-4} \cong
H^{1}(X_{o},{\mathbb W}_{u})$. Note that ${\mathfrak g}$ is
semisimple, by general theory.  In Lemma~\ref{lem-irreducible} we
proved that $H^{1}(X_{o},{\mathbb W}_{u})$ is an irreducible
representation of ${\mathfrak g}$. The above proof works for the Lie
algebra since the $D_{ij}^2$ are unipotent matrices with Jordan blocks
of length at most one and thus $A_{ij}:= D_{ij}^2 -1$ are elements of
${\mathfrak g}$.  Hence the above proof shows that the $A_{ij}$ act
irreducibly.

For some $i,j$ we have $A_{ij}$ decomposing into two Jordan blocks of
length one (this is the case for $A_{1,2}$, $A_{3,4}$ etc.). However
there are some $i,j$ where the monodromy $M_{i,j}$ is a reflection
(for example $i = 2, j =11$), where the $A_{i,j}$ has a single Jordan block
of length one.

Next, note that by isolating the monodromy representation on the part
of the curve where the local system is trivial, we obtain a monodromy
representation of the direct sum of two copies of the cohomology of a
hyperelliptic curve of genus $g'$ with the monodromy acting
diagonally.  If we take all branch points for $i\geq 11$ then this has
genus $g' = g-5$. Deligne's argument from
\cite[Section~4.4]{deligne-weil2} or the argument from \cite{janssen},
works for the hyperelliptic monodromy action on the standard
cohomology (the monodromy is generated by conjugate Dehn twists) so
this monodromy group is $\op{Sp}(2g-10)$.  In particular we have
\[
\liesp(2g - 10) \subset {\mathfrak g}\subset \liesp(4g-4)
\]
where the composite inclusion is the linear embedding of the diagonal
action on the direct sum of two copies of the standard representation
of $\liesp(2g - 10)$.

This shows that our monodromy action satisfies the hypothesis of the
following purely algebraic theorem:

\begin{theo} \label{theo-nested} 
There exists $g_0$ with the following property:
suppose $g \geq g_0$ and suppose ${\mathfrak g}$ is a semisimple Lie algebra
sitting in a pair of inclusions
\[
\liesp(2g - 10) \subset {\mathfrak g}\subset \liesp(4g-4)
\]
where the composite inclusion is the linear embedding of the diagonal
action on the direct sum of two copies of the standard representation
of $\liesp(2g - 10)$. Suppose that the action of ${\mathfrak g}$ on
${\mathbb C}^{4g-4}$ is irreducible, and suppose ${\mathfrak g}$ contains an
element $A$ 
which acts on ${\mathbb C}^{4g-4}$ with a single Jordan block of length
one. Then ${\mathfrak g}= \liesp(4g-4)$. 
\end{theo}
{\bf Proof.} We first claim that ${\mathfrak g}$ is simple. If that
were not the case, then we could write ${\mathfrak g} = {\mathfrak
g}_1 \times {\mathfrak g}_2$ and so the representation $W={\mathbb
C}^{4g-4}$ would decompose as an exterior tensor product
$=W_1\boxtimes W_2$ of representations of ${\mathfrak g}_1$ and
${\mathfrak g}_2$. This however can be ruled out by looking at the
element $A$.  If it is nontrivial in both factors then it would act on
$W$ by a Jordan normal form which is the tensor product of two
nontrivial Jordan normal forms, in particular it would have a Jordan
block of length $>1$; if it was nontrivial in only one of the factors
then it would act by the tensor product of a nontrivial Jordan form,
by a trivial vector space (of dimension $>1$); thus it would have at
least two Jordan blocks. In either case this contradicts the
hypothesis that $A$ acts on ${\mathbb C}^{4g -4}$ with a single Jordan
block of length one. This proves that ${\mathfrak g}$ is simple.

Now we can choose $g_0$ big enough so that the dimension of $\liesp(2g
- 10)$ is bigger than the dimension of the exceptional simple Lie
algebras.  By classification, this means that ${\mathfrak g}$ is of
one of $\lieso(2m)$, $\lieso(2m -1)$, $\liesp(2m)$, or
${\liesl}(m)$. Looking at dimensions, we get $m \geq g - C$ where $C$
is some constant.

On the other hand, note that all the fundamental representations of
the classical groups are essentially (the only exception being the
spin representation) the wedge powers of the standard representation.
Combined with Weyl's dimension formula this implies that there is an
$m_0$ such that for $m\geq m_0$ the only irreducible representations
of dimension $<5m$ of one of the classical groups above, are the
fundamental representation or (in the last case) the dual of the
fundamental representation.

This claim gives that ${\mathfrak g}$ acts by the standard representation,
which immediately implies that it is equal to $\liesp(4g-4)$ (it
can't be orthogonal or special linear because we already know it is
contained in the symplectic group). \eop

\

\bigskip

\noindent
We are now in a position to complete the

\

\medskip

\noindent
{\bf Proof of Theorem~\ref{thm-finite-index}:} Let $g_{0}$ be such
that Theorem~\ref{theo-nested} holds. In view of
Corollary~\ref{cor-localization} and Lemma~\ref{lem-condition-ii} we
only need to show that the Zariski closure $G$ of the monodromy action
of $\pi_{1}(B,o)$ on $H_{1}(X_{o},\op{ad}(\rho))$ acts with an open
orbit on $H_{1}(X_{o},\op{ad}(\rho))$.  

By Theorem~\ref{theo-nested}, we
get that each subspace $H_{1}(X_{o},{\mathbb W}_u) \subset
H_{1}(X_{o},\op{ad}(\rho))$ yields a monodromy representation equal to
the symplectic group.  Note that Theorem~\ref{theo-nested}
implies that the
Lie algebra of the monodromy group is equal to ${\liesp}(4g-4)$. However,
in view of the fact that we already have an inclusion of the monodromy
group of $H_{1}(X_{o},{\mathbb W}_u)$ in $\op{Sp}(4g-4)$, we get that this
inclusion is surjective so the monodromy group is equal to $\op{Sp}(4g-4)$
which we now write as $\op{Sp}(H_{1}(X_{o},{\mathbb W}_u))$.

In particular we have a natural
inclusion 
\[
G \subset \op{Sp}(H_{1}(X_{o},{\mathbb C})) \times \prod_{u\neq (\bone,0)}
\op{Sp}(H_{1}(X_{o},{\mathbb W}_u)),
\]
so that the projection on each factor is surjective.  On the other
hand, $G$ is semisimple.  Going back to the level of Lie algebras,
this implies that the simple summands of ${\mathfrak g}$ are
$\liesp(2g)$ occuring once, and $\liesp(H_{1}(X_{o},{\mathbb W}_u)) =
\liesp(4g-4)$ occuring a certain number of times. Call these summands
${\mathfrak s}_1,\ldots , {\mathfrak s}_k$.

Each irreducible factor $H_{1}(X_{o},{\mathbb W}_u)$ in the
representation $H_{1}(X_{o},\op{ad}(\rho))$, is an irreducible 
representation of the Lie algebra
\[
\liesp(2g) \oplus {\mathfrak s}_1\oplus \ldots \oplus {\mathfrak s}_k.
\]
As such, it decomposes {\em a priori} into an exterior tensor product
of representations of the summands; but by dimension considerations,
this tensor product must just be an irreducible representation of one
of the summands ${\mathfrak s}_i$.  Also this representation is isomorphic to
the standard representation of ${\mathfrak s}_i = \liesp(4g-4)$.  Thus
each $H_{1}(X_{o},{\mathbb W}_u)$ comes from the standard
representation composed with a 
projection onto one of the factors.  Note also that the representation
${\mathbb C}^{2g}$ comes from the standard representation composed with the
projection onto the factor $\liesp(2g)$.

The above statements on the level of Lie algebras imply the same
things for the connected components of the Lie groups. We obtain that
the connected component of the monodromy group $G$ decomposes as a
product
\[
G^o = \op{Sp}(2g) \times S_1 \times \ldots \times S_k
\]
where each $S_i$ is equal to $\op{Sp}(4g-4)$ and $S_i$ acts on a direct sum
of $r_i$ copies of its standard representation.

Our goal now is to prove that all of the $r_i$ are equal to
$1$. Suppose the contrary, i.e. suppose that there are two distinct
components $u$ and $u'$ such that the same $S_i$ acts on
$H_{1}(X_{o},{\mathbb W}_u)$ and 
$H_{1}(X_{o},{\mathbb W}_{u'})$. In particular this means that
$H_{1}(X_{o},{\mathbb W}_u)$ 
and $H_{1}(X_{o},{\mathbb W}_{u'})$ are 
isomorphic as representations of $G^o$.

The elements $D_{ij}^2$ are unipotent so they go into $G^o$.  Let
$\Gamma$ denote the subgroup of the monodromy group which maps into
$G^o$.  We obtain a map
\[
{\mathbb C} [\Gamma ] \rightarrow {\mathbb C} [ G^o].
\]
In particular the action of the group algebra ${\mathbb C} [\Gamma ]$
on $\bigoplus_u H_{1}(X_{o},{\mathbb W}_u)$ factors through the action
of the group algebra of $G^o$.  Thus, with our assumption of the
previous paragraph that some $r_i$ is $>1$, we would get two
components $H_{1}(X_{o},{\mathbb W}_u)$ and $H_{1}(X_{o},{\mathbb
W}_{u'})$ which are isomorphic as representations of the group algebra
${\mathbb C} [\Gamma ]$.  Therefore for every element $E\in {\mathbb
C} [\Gamma ]$, the eigenvalues of $E$ acting on $H_{1}(X_{o},{\mathbb
W}_u)$ and $H_{1}(X_{o},{\mathbb W}_{u'})$ are the same. We will write
down elements $E$ which act with single nonzero eigenvalues; thus
these values are the same for $u$ and $u'$. We will show that this
implies that $u$ and $u'$ are the same component.  This will
contradict the assumption that $r_i > 1$ so it will prove the
statement that all of the $r_i$ are equal to $1$.

We will write down our $E$ as products of elements of the form
$A_{ij}:=D_{ij}^2-1$.  The component $u$ (respectively $u'$) is
determined by the numbers $b$ and $c$ (respectively $b'$ and $c'$)
which occur above. These are taken modulo the order $n$ of the root of
unity $\gamma$, and interchanging $b\leftrightarrow -b$ and
$c\leftrightarrow -c$ doesn't change the dihedral component $u$. Thus,
proving that the two components are the same means that we want to
show
\[
(b,c) = \pm (b',c')\;\; \mbox{in} \;\; ({\mathbb Z} /n)^2.
\]

Let $i,j=7,9,11$. 
Look for example at 
\[
A_{2,i}A_{1,2}.
\]
This takes the vector $v_{2,j}L_{2,j}$ first to $v_{2,j}L_{1,2}$ and
then to
\[
(1 + P^{-x_i}R) v_{2,j}L_{2,i}
\]
where $x_i=b,c$ or $0$ depending on whether $i=7,9$ or $11$.  
Look now at 
\[
E = A_{2,i}A_{1,2}A_{2,j}A_{1,2}.
\]
It has an image vector which is a multiple of $v_{2,i}L_{2,i}$ so this
can be its only nonzero eigenvector.  Its action on this vector is
\[
Ev_{2,i}L_{2,i} = (1+P^{-x_i}R)(1 + P^{-x_j}R)v_{2,i}L_{2,i}.
\]
Note also that the matrix $(1+P^{-x_i}R)$ is itself of rank one so the
product of matrices appearing above also has a single nonzero
eigenvalue.  In particular the unique nonzero eigenvalue of $B$ is
equal to the unique nonzero eigenvalue of the matrix $(1+P^{-x_i}R)(1
+ P^{-x_j}R)$.  This matrix may be written as
\[
\left( \begin{array}{cc} 
1 + \gamma ^{x_j-x_i} & \gamma ^{-x_i} + \gamma ^{-x_j} \\
\gamma ^{x_i} + \gamma ^{x_j} & 1 + \gamma ^{x_i-x_j}   
\end{array} \right) .
\]
Its eigenvector is 
\[
v_{2,i} = (1+P^{-x_i}R)v_+ = \left( \begin{array}{c} 
1   \\
\gamma ^{x_i}  
\end{array} \right) .
\]
Calculating
\[
\left( \begin{array}{cc} 
1 + \gamma ^{x_j-x_i} & \gamma ^{-x_i} + \gamma ^{-x_j} \\
\gamma ^{x_i} + \gamma ^{x_j} & 1 + \gamma ^{x_i-x_j}   
\end{array} \right) 
\left( \begin{array}{c} 
1   \\
\gamma ^{x_i}  
\end{array} \right) 
=
\left( \begin{array}{c} 
2 + \gamma ^{x_j-x_i} + \gamma ^{x_i - x_j}   \\
2\gamma ^{x_i} + \gamma ^{x_j} + \gamma ^{2x_i -x_j}
\end{array} \right) 
\]
says that the eigenvalue is equal to 
\[
2 +  \gamma ^{x_j-x_i} + \gamma ^{x_i - x_j}   .
\]  

Now take various values of $i$ and $j$, and compare the results for the
components $u$ and $u'$. We obtain:
\begin{description}
\item[from $i=11,j=7$,]
\[
\gamma ^b + \gamma ^{-b} = \gamma ^{b'} + \gamma ^{-b'};
\]
\item[from $i=11,j=9$,]
\[
\gamma ^c + \gamma ^{-c} = \gamma ^{c'} + \gamma ^{-c'};
\]
\item[and from $i=9,j=7$,] 
\[
\gamma ^{b-c} + \gamma ^{c-b} = \gamma ^{b'-c'} + \gamma ^{c'-b'}.
\]
\end{description}

In ${\mathbb Z} /n$ these equations give:
\[
b = \pm b',\;\; c = \pm c', \;\; (b-c) = \pm (b' - c').
\]
The first two equations admit four possibilities: either 
\[
(b,c) = (b', c'), \;\; \mbox{or} \; (b,c) = -(b',c'),
\]
or else
\[
(b,c) = \pm (b', -c').
\]
The first two possibilities are what we want to show. In the last two
possibilities we have
\[
b-c = \pm (b' + c').
\]
Thus the third equation above says
\[
b' + c' = \pm (b' - c').
\]
This says that either $c' = -c'$ or else $b' = -b'$. In either case, the
equation $(b,c) = \pm (b', -c')$ gets transformed into the equation
$(b,c) = \pm (b' , c')$ so we are done.  In fact, we could have noted
here  that since $n$ is odd, the equations $c'=-c'$
or $b'=-b'$ do not occur, but the present argument works even when $n$
is even.

We have now shown that if $H_{1}(X_{o},{\mathbb W}_u)$ and $H_{1}(X_{o},{\mathbb W}_{u'})$ are isomorphic as
representations of $G^o$ then $u$ and $u'$ represent the same dihedral
component of our representation.  This implies that all of the $r_i$
are equal to one, which in turn gives that
\[
G^o = \op{Sp}(2g) \times \prod_{u \neq (\bone,0)} \op{Sp}(H_{1}(X_{o},{\mathbb W}_u)).
\]
(Actually, the same is true of the full monodromy group because the
full group is also contained in this product of symplectic groups:
$G=G^o$.)  This group acts on its representation ${\mathbb C} ^{2g}
\oplus \bigoplus_{u \neq (\bone,0)} H_{1}(X_{o},{\mathbb W}_u)$ with
an open orbit. Combined with Lemma~\ref{lem-condition-ii} and
Theorem~\ref{theo-alternative} this completes the proof of
Theorem~\ref{thm-finite-index}. \eop

\subsection{Lefschetz pencils} \label{sec-lefschetz} 

We would now like to extend the techniques of the previous section in
order to prove Theorem~\ref{thm-lefschetz-pencil}.

Let $Z$ be a smooth projective surface with $b_{1}(Z) = 0$. Let
${\mathcal O}_{Z}(1)$ be a very ample line bundle on $Z$ and let
${\mathbb P}^{1} \subset {\mathbb P}(H^{0}(Z,{\mathcal O}_{Z}(k)))$ be
a generic line. Denote by $\varepsilon : \widehat{Z} \to Z$ the
blow-up of $Z$ at the base points of the pencil of curves $\{ D_{t}
\}_{t \in {\mathbb P}^{1}}$, and let $f : \widehat{Z} \to {\mathbb
P}^{1}$ be the corresponding Lefschetz fibration. Let $p_{1}, \ldots,
p_{\mu} \in {\mathbb P}^{1}$ be the critical points of $f$ and let $B
= {\mathbb P}^{1}-\{p_{1}, \ldots, p_{\mu}\}$ and $X = f^{-1}(B)$. Let
$X_{t}$, $t \in B$ denote the fiber of $f$ over $t$, or equivalently,
the strict transform of the divisor $D_{t} \subset Z$.

Fix $g_0$ so that Theorem~\ref{theo-nested} applies. Then, as we saw
at the end of the previous section, it follows that the monodromy
group of the hyperelliptic family of genus $g_0$ acting on the
cohomology of the full local system $\op{End}({\mathbb V}_{n})$
corresponding to $\op{ad}(\rho)$ is equal to \linebreak 
$\op{Sp}(2g_0) \times \prod _u
\op{Sp}(H_{1}(X_{o},{\mathbb W}_{u}))$ with $\op{Sp}(H_{1}(X_{o},{\mathbb
W}_{u}))=\op{Sp}(4g_0-4)$.  We also assume that $g_0 \geq 6$, so there are
several branch cuts along which the representation $\rho$ is
the identity.  Set $m:= 2g_0+1$.

Let ${\mathcal L} := {\mathcal O}_{Z}(k)$ denote our line bundle. Our
assertions will be made for $k$ big enough.  Let $\bbD \subset
{\mathbb P} H^0(Z,{\mathcal L} )$ denote the discriminant locus
consisting of the sections defining singular curves.  We will fix a
base point $o \in {\mathbb P} H^0(Z,{\mathcal L} )\!-\!\bbD$ (chosen
specially below); and as always $X_o$ will denote the smooth curve
defined by the section $o$. Then $ \pi_1({\mathbb P}(H^0)\!-\!\bbD,
o)$ acts by diffeomorphisms on $X_o$ and hence it acts on
$M_B(X_o,n)$. Furthermore, by the Lefschetz hyperplane section
theorem, the geometric monodromy action $\pi_{1}(B,o) \to
\op{Map}(X_{o})$ for the family $f : X \to B$ factors through the
natural map $\pi_{1}(B,o) \to \pi_1({\mathbb P}(H^0)\!-\!\bbD, o)$ and
so it suffices to show that $\pi_1({\mathbb P}(H^0)\!-\!\bbD, o)$ acts
on $M_B(X_o,n)$ with a Zariski dense orbit.

As before we shall fix a local system $\rho$ on $X_{o}$, with finite
monodromy factoring through a representation of the dihedral
Heisenberg group. Then there is a subgroup of finite index in
$\pi_1({\mathbb P}H^0\!-\!\bbD,o)$ which preserves $\rho$, so it
acts on the space
\[
T_{[\rho]} M_{B}(X_{o},n) = H^1(X_o,\op{ad}(\rho))= H^1(X_o,{\mathbb
C} )\oplus \left( \bigoplus_{u \neq (\bone,0)} H^1(X_o, {\mathbb
W}_u)\right).
\]
According to Corollary~\ref{cor-localization} and
Lemma~\ref{lem-condition-ii} we only need to show that the Zariski
closure $G$ of the image of this monodromy action acts with an open
orbit on $H^1(X_o,\op{ad}(\rho))$. Our technique will be to show that
$G$ is as big as possible, given the above decomposition and the fact
that it preserves symplectic forms on everything. To achieve this we
will use a family of curves in the linear system $|{\mathcal L}|$ which
have hyperelliptic handles and will apply to the curves the results
for the hyperelliptic case obtained in
Section~\ref{sec-hyperelliptic}.

We will define a particular subspace $\bbE \subset {\mathbb P}
H^0(Z,{\mathcal L} )$, and among other things choose $o\in \bbE\!-\!\bbD
\cap \bbE$.  Then the fundamental group $\pi_1(\bbE\!-\!\bbD \cap \bbE,o)$ is
contained in $\pi _1(H^0\!-\!\bbD,o)$ and again a finite index
subgroup will act on $H^1(X_{o},\op{ad}(\rho))$. The subspace $\bbE$ will be
designed so that $\pi_{1}(\bbE\!-\!\bbD \cap \bbE,o)$ preserves the handle
decomposition of $X_{o}$. Using this we will obtain first a smaller
subgroup of $G$ and then apply an argument using loops in the full space
${\mathbb P} H^0(Z,{\mathcal L} )-\bbD$ (which don't preserve the
handle decomposition) to obtain the full group $G$.

Fix a point $P\in Z$.
Choose a collection of sections $a_{0}, \{f_{i}\}_{i=1}^{m},
\{s_{j}\}_{j=0}^{N} \in H^{0}(Z,{\mathcal L})$ so that for a suitable
local trivialization of ${\mathcal L}$ and local
coordinates $(x,y)$ near $P$ we have:

\begin{description}
\item[(a)] near $P$ we have $a_0=x^m -y^2$;
\item[(b)] the section $a_0$ has no singularities other than $P$;
\item[(c)] near $P$ the sections $f_i$ have the form $x^i$;
\item[(d)] the sections $s_{j}$ vanish to order at least $m+2$ at $P$
and together with $a_0$ form a basis for a linear 
system which has no base points outside of $P$.
\end{description}

\

\noindent
It is clear that by taking $k$ big enough we can always find such
$a_{0}$, $f_{i}$ and $s_{j}$. In addition we will need to choose the
$s_{j}$'s so that they satisfy certain connectedness conditions which
we shall describe further on (Lemma \ref{connectednessI} and Lemma
\ref{interchange}).

With these choices we let $\bbE$ be the affine space of sections of
${\mathcal L}$ of the form 
\[
a_0 + \sum_{i=0}^{m-1} t_i f_i + \sum _{j=0}^N z_j s_{j}, 
\]
and we take as a base point the point $o \in \bbE$ corresponding to
the values $z_j=0$ and $t_i=0$ for $i>0$, with $t_0=1$.  Write this as
$o = (1,0,\ldots , 0)$.  We will work on open polydisks in $\bbE$ of
the form $|t_{0} - 1| < \bB$, $|t_i| < \bB$ for $i > 0$ and $|z_j| < \bC$
for all $j$ (with $\bB$ and $\bC$ to be determined later), possibly after
rescaling the $f_i$.

To achieve the desired behavior of the monodromy we start by looking
at the choice of the $s_{j}$ and $\bC$.  For $k$ big enough we can
choose a linearly independent family of sections $s_{j}$ such that the
$s_{j}$ all vanish to order $>m$ at $P$, and such that the linear
system they generate is without base points away from $P$.  Let
$\bbG\subset {\mathbb P} H^0(X,{\mathcal L} )$ denote the subspace
generated by the $s_{j}$ and by $a_0$.  It is a projective space, with
a codimension one projective subspace $\bbG_{\infty} \subset \bbG$
corresponding to the linear system spanned by the $s_{j}$ without
$a_0$.  The family of sections $a_0 + \sum_j z_j s_{j}$ provides a
system of affine coordinates $z_j$ for the complementary affine space
${\mathbb A}_\bbG := \bbG-\bbG_{\infty}$.  On the other hand, over
${\mathbb A}_\bbG$ the universal family of curves is a family which is
holomorphically locally trivial near the singular point $P$. Indeed,
over any small enough disc in the coordinates $z_j$, one can choose
local coordinates at $P$, depending on the $z_j$'s, such that $a_0 +
\sum_j z_j s_{j}$ has the form $x(z_j)^m -y(z_j)^2$. Let $\bbD_\bbG$
denote the subset of points in $\bbG$ parametrizing curves that have
singularities outside of $P$, union the $\bbG_{\infty}$ which corresponds
to curves with bigger singularities than usual at $P$. Over $\bbG-\bbD_\bbG$
we obtain a family of curves which are smooth except for their
singularities at $P$, and the family is holomorphically locally
trivial (hence topologically locally trivial) along the section
corresponding to the point $P$.

Let $s'$ be a point in the complement $\bbG-\bbD_\bbG$ and let
$X_{s'}$ be the fiber over $s'$.  With $m$ odd, the singularity of
each fiber at $P$ is a higher-order cusp, so in fact the fibers such
as $X_{s'}$ are singular curves which are topologically (but not
differentially) manifolds.  The local topological triviality of the
family means that it makes sense to speak of the monodromy action of
$\pi _1(\bbG-\bbD_\bbG, s')$ on the cohomology of the fiber $X_{s'}$
with coefficients in the trivial local system.

Now we come to the first of the connectedness conditions referred to
above.  In fact the subject of the monodromy of pencils having
singularities especially at the base locus, has been intensively
studied recently notably by Tibar \cite{Tibar1,Tibar2}. Our
situation above is a very special easy case of this phenomenon so we
don't need to call upon his general results.

\begin{lem}
\label{connectednessI}
For $k$ big enough and by choosing a big enough family of sections
$s_{j}$, we can ensure that the Zariski closure of the monodromy
action of $\pi_1(\bbG-\bbD_\bbG,s')$ on $H^1(X_{s'}, {\mathbb C} )$ is the
full symplectic group.
\end{lem}
{\bf Proof.}  If we choose a general line ${\mathbb A}^1\subset
{\mathbb A}_\bbG$ then by standard Lefschetz theory (see
e.g. \cite{katz-sga}, \cite{looijenga-notes}) the fundamental group of
$\bbG-\bbD_\bbG$ is generated by the loops in ${\mathbb A}^1 -
{\mathbb A}^1\cap \bbD_\bbG$.  Here (and this is the important point
of the argument) we can take only the loops which go around the points
in ${\mathbb A}^1\cap \bbD_\bbG$; we don't need to look at the loop going
around the point at infinity since it is the product of the others.
Thus the fundamental group of $\bbG-\bbD_\bbG$ is generated by loops going
around the affine part of the discriminant $\bbD_\bbG -\bbG_{\infty} =
\bbD_\bbG \cap {\mathbb A}_\bbG$.

For $k$ big enough and by choosing a big enough family of sections
$s_{j}$, the affine part of the discriminant divisor $\bbD_\bbG
-\bbG_{\infty}$ is irreducible, hence connected.  Note that we could
never assure that $\bbD_\bbG$ is connected since it contains
$\bbG_{\infty}$ as an irreducible component --- thus the importance of
saying that the monodromy is generated by loops around the affine
piece.  To get this connectedness we follow the standard argument in
the theory of Lefschetz pencils: the discriminant divisor in the
affine piece is the image of an affine space bundle over the surface
$Z-P$.  For big enough values of $k$ this family of affine spaces
(which to a point $x\in Z-P$ associates the subspace of sections in
${\mathbb A}_\bbG$ which are singular at $x$) is a vector bundle over
$Z-P$; thus its image is irreducible.  Of course we also choose $k$ so
that the general point in this divisor corresponds to an ordinary
double point of the curve.

Now the Kazhdan-Margulis result as reported by Deligne
\cite[5.10]{deligne-weil2} works the same way to show that the
monodromy of the fundamental group of $\bbG-\bbD_\bbG$ on the cohomology
of the family of curves with trivial coefficients, has Zariski closure
equal to the full symplectic group.  Indeed the monodromy around a
point of $\bbD_\bbG-\bbG_{\infty}$ is a symplectic transvection (because
the singular curve has an ordinary node), and the connectedness of the
divisor means that all of these elements are conjugate. As we have
seen above, they generate the monodromy group, so we have a group
generated by a family of conjugate symplectic transvections.
Furthermore the monodromy representation has no fixed vectors (a fixed
vector would correspond to a class in $H^1(Z,{\mathbb C})$ which we
have assumed is trivial). \eop

\bigskip

\noindent
Choose $k$ and the $s_{j}$ as per the above lemma.  Choose an explicit
collection of loops $\gamma _k$ in $\bbG-\bbD_\bbG$ which generate the
monodromy, and choose $\bC$ big enough so that these loops are contained
in the region $|z_j| < \bC$.

Our coordinate patch around $P$ will consist of a nested pair of balls
$U \subset U' \subset Z$ together with a pair of coordinate functions
$(x,y): U' \rightarrow {\mathbb C}^2$ sending $U$ (respectively $U'$) to the
ball of radius $\bT$ (respectively $\bT'$) in ${\mathbb C}^2$. We can assume
that $x$ and $y$ come from sections of ${\mathcal O}_Z(k_0)$ for some
$k_0$, via a trivialization of this last line bundle over
$U'$. Furthermore the trivialization can be assumed to come from a
given section $u$ of ${\mathcal O}_Z(1)$ which doesn't vanish at $P$.
Then provisionally put
\[
f_i := x^iu^{k-ik_0}\;\;\; \mbox{for}\;\;\; 0\leq i \leq m-1 \;\;\;
\mbox{and} \;\;\; a_0 := x^m u^{k - mk_0} - y^2 u^{k-2k_0}.
\]
These come from global sections of ${\mathcal O}_{Z} (k)$ (choose $k >
mk_0$).

\begin{lem}
\label{rescaling}
For any fixed choice of the sections $s_{j}$ and $a_0$, of the
constants  $\bT$ and $\bT'$, and for
any $\delta >0$, we can make a rescaling of our local coordinates and
of the $f_i$ so that the following hold: 
\begin{itemize}
\item we can retain the properties
(a) and (c), while $a_0$ and the $s_{j}$ remain fixed in
$H^0(Z,{\mathcal L} )$; 
\item the sections $f_i$ become arbitrarily small
inside $H^0(Z,{\mathcal L} )$; 
\item the coordinate patches $U$ and $U'$
become arbitrarily small inside $Z$, and
\[
\left| \sum _j z_j s_{j} \right| < \delta 
\]
on the coordinate patch $U'$, for all $|z_j|<\bC$.
\end{itemize}
\end{lem}
{\bf Proof.}  We can rescale $x$ by a factor of $\lambda^2$ and $y$ by
a factor of $\lambda^m$, and the trivialization by a factor of
$\lambda^{-2m}$ (in other words scale $u$ by a factor of
$\lambda^{2m/k_0}$).  We retain the expressions (a) and (c), $a_0$
remains fixed in $H^0(Z,{\mathcal L} )$, and the sections $f_i$ become
arbitrarily small inside $H^0(Z,{\mathcal L} )$. Also the coordinate
patches $U$ and $U'$ become arbitrarily small inside $Z$. Finally,
since the $s_{j}$ (which are fixed) vanish to order $>m$ at $P$, for
any $\delta >0$ we can choose the rescaling so that the required
estimate holds. \eop

\

\bigskip

\noindent
The family of curves of the form 
\[
a_0 + \sum _{i=0}^ {m-1} t_if_i = x^m - y^2 + \sum _{i=0}^ {m-1} t_ix^i
\]
gives the full family of hyperelliptic curves in our coordinate patch.
In particular the monodromy of this family (i.e. the fundamental group
of the complement of the discriminant locus) acts as the braid group
of braids on $m$ strands.  Choose loops generating
this monodromy and choose $\bB$ so that the loops are contained in the
region $|t_i|< \bB$.

Let $\partial U$ be the spherical boundary of the coordinate patch.
In terms of the local coordinates the equation for $\partial U$ is
$|(x,y)| = \bT$.  For $\bT$ big enough, the intersection of the curve
\[
a_0 + \sum _{i=0}^ {m-1} t_if_i
\]
with the sphere $\partial U$ will remain approximately equal to a
single fixed circle as $t_i$ vary in the region $|t_i|<\bB$.  As pointed
out above, if we made a sufficient rescaling at the start then the
ball $U$ will still be small in $Z$.  Furthermore, by the estimate of
Lemma \ref{rescaling}, addition of the terms $\sum _j z_j s_{j}$ will
not move the intersection circle by very much either.

We have thus found parameters for our family of curves $\bbE$ (and a
rescaling of the $f_i$) such that for any point $(t,z)$ in the family
satisfying the bounds $|t_i|<\bB$ and $|z_j|<\bC$, the intersection
$X_{(t,z)}\cap \partial U$ remains close to a single fixed circle.

Finally we choose as a base point $o=(1,0,\ldots , 0)$ i.e. $t_0=1$.
By choosing $\delta$ small enough in the above choices, we can insure
using the estimate of Lemma \ref{rescaling} that when we let the $z$
coordinates go around the loops $\gamma_k$, the piece of the curve
inside $U$ doesn't move too far from the base curve and in particular
the monodromy action on this piece is trivial.  Similarly, note that
with a very small scaling factor $\lambda$ the sections $f_i$ become
very small compared to $a_0$ (or more precisely, they become small
compared to the differential of $a_0$ along the zero-set of $a_0$).
So when $t$ goes around loops generating the braid group action, this
doesn't move very much the curve outside of $U'$. Furthermore we choose
$U$ so that these loops act trivially in the coordinate patch on
$U'-U$. Thus when $t$ moves we obtain a braid group action on the
piece of the curve inside $U$ and a trivial action on the piece
outside $U$.  We can sum all this up as follows:

\begin{cor}
\label{sumup}
Over $\bbE$ the family of curves decomposes into a family of
hyperelliptic curves of genus $g_0$ joined onto a family of curves of
genus $g-g_0$ along a circle which stays essentially fixed.  There is
a collection of paths in $\bbE$ which generate a monodromy action on
the genus $g-g_0$ piece whose Zariski closure is the full symplectic
group of the cohomology of the genus $g-g_0$ curves (Lemma
\ref{connectednessI}).  These paths act trivially on the hyperelliptic
piece. On the other hand, there are paths in $\bbE$  generating the braid
group action on the hyperelliptic piece, which in turn act trivially
on the piece of genus $g-g_0$.
\end{cor}

\

\bigskip

\noindent
Next, fix our base representation $\rho$ to be trivial on the genus
$g-g_0$ piece and equal to the dihedral Heisenberg representation
chosen in the hyperelliptic argument of the previous sections for the
hyperelliptic genus $g_0$ piece.  Note that since $m=2g_0+1$ rather
than $2g_0+2$, one of the branch points on the hyperelliptic handle is
at infinity; we assume that this branch point is part of a branch cut
on which the passing matrix is trivial. Furthermore there is at least
one other branch cut on which the passing matrix is trivial too.

If we consider monodromy elements coming from the family $\bbE$, these
preserve the cutting-up of our curve into pieces $X_{o}\cap U$ of genus
$g_0$ and $X_{o} -X_{o}\cap U$ of genus $g-g_0$. This monodromy group
preserves the decomposition of the cohomology of $\op{ad}(\rho )$ into a
sum of two pieces, one of dimension $2n^{2}(g-g_0)$ and the other of
dimension
\[
2(n^2-1)(g_0 -1) + 2g_0. 
\]
The first piece corresponds to the monodromy on the sum of $n^2$
copies of the trivial representation, since $\rho$ and thus
$\op{ad}(\rho )$ are trivial outside $X_{o}\cap U$.  Now the
cohomology of $X_{o} -X_{o}\cap U$ with coefficients in the trivial
local system is isomorphic to the cohomology of the Riemann surface
obtained by glueing in a disc along the boundary $X_{o}\cap \partial
U$. In turn, this Riemann surface is homeomorphic to the cusp curve
considered in the family parametrized by ${\mathbb G}$ above. Thus the
monodromy result of Lemma \ref{connectednessI} implies that the
monodromy action of the loops $\gamma_k$ on the cohomology of $X_{o}
-X_{o}\cap U$ with coefficients in the trivial local system, has
Zariski closure equal to the symplectic group $\op{Sp}(2(g-g_0))$.  Now
taking into account the fact that the restriction of ${\mathbb V}_{n}$ to
$X_{o} -X_{o}\cap U$ is the direct sum of $n^2$ copies of the trivial
representation; we get that the monodromy action on this first piece
is a diagonal copy of $\op{Sp}(2(g-g_0))$ embedded in the product of $n^2$
copies of its standard representation.

The second piece corresponds to the monodromy action on the cohomology
of our genus $g_0$ handle with coefficients in $\op{ad}(\rho )$.  We have
an action of the braid group $B_m$ on this second piece equal to the
braid monodromy action for the family of hyperelliptic curves.  Of
course not all elements of $B_m$ preserve the representation $\rho$,
and as before we look only at elements which preserve $\rho$. The
results of Section~\ref{sec-hyperelliptic}  apply here, giving the Zariski
closure of the monodromy of this braid action on the cohomology of the
hyperelliptic piece.

We are now in the following situation. Consider the cohomology
$H^{1}(X_{o},\op{ad}(\rho))$. Let $H^{1}(X_{o},\op{ad}(\rho))^{\rm
hyper}$ and $H^{1}(X_{o},\op{ad}(\rho))^{g-g_0}$ respectively denote
the cohomologies of $\op{ad}(\rho)$ over the hyperelliptic piece and
the complementary piece. Note that in $H^{1}(X_{o},\op{ad}(\rho))^{\rm
hyper}$ we restrict to classes which are zero on the boundary circle
joining the two pieces.  By Mayer-Vietoris, restriction gives an isomorphism 
\[
H^{1}(X_{o},\op{ad}(\rho)) \stackrel{\cong}{\rightarrow}
H^{1}(X_{o},\op{ad}(\rho))^{\rm hyper} \oplus
H^{1}(X_{o},\op{ad}(\rho))^{g-g_0}.
\]
Furthermore both pieces have natural symplectic forms and the
symplectic form on \linebreak  $H^{1}(X_{o},\op{ad}(\rho))$ corresponds to the
direct sum. 

Combined with the decomposition 
\[
H^{1}(X_{o},\op{ad}(\rho)) = H^{1}(X_{o},{\mathbb C})\oplus \left(
\bigoplus_{u \neq (\bone,0)} H^{1}(X_{o},{\mathbb W}_{u}) \right)
\]
this yields splittings
\[
H^{1}(X_{o},\op{ad}(\rho))^{\rm hyper} = H^{1}(X_{o},{\mathbb C})^{\rm
hyper} \oplus \left( \bigoplus_{u \neq (\bone,0)} H^{1}(X_{o},{\mathbb
W}_{u})^{\rm hyper} \right),
\]
and similarly
\[
H^{1}(X_{o},\op{ad}(\rho))^{g-g_0} = H^{1}(X_{o},{\mathbb C})^{g-g_0}
\oplus \left( \bigoplus_{u \neq (\bone,0)} H^{1}(X_{o},{\mathbb
W}_{u})^{g-g_0} \right)
\]
Again these decompositions are compatible with the symplectic form.
The action of \linebreak $\pi_1(\bbE -\bbE\cap \bbD, o)$ on
$H^{1}(X_{o},\op{ad}(\rho))$ preserves this decomposition.

Let $G$ denote the global monodromy group i.e. the complex Zariski
closure of the monodromy image of $\pi_1({\mathbb P} H^0(Z,{\mathcal
L} ) - \bbD,o)$ acting on $H^{1}(X_{o},\op{ad}(\rho))$. Let $G^{\rm
hyper}$ (respectively $G^{g-g_0}$) denote the Zariski closures of the
monodromy of the family $\bbE$ (i.e. of the image of $\pi_1(\bbE
-\bbE\cap \bbD,o)$) acting on each of the pieces in the above
decomposition.

\begin{lem}
\label{product}
With these notations,
the product group is contained in the global monodromy:
\[
G^{\rm hyper}\times G^{g-g_0} \subset G. 
\]
Furthermore,
\[
G^{\rm hyper} = \op{Sp}(2g_0) \times \prod_{u \neq (\bone,0)}
\op{Sp}(H^{1}(X_{o},{\mathbb 
W}_{u})^{\rm hyper}),
\]
and $G^{g-g_0}$ is the diagonal copy of $\op{Sp}(2(g-g_0))$ acting on 
$H^{1}(X_{o},{\mathbb C})^{g-g_0}\cong {\mathbb C} ^{2n^2(g-g_0)}$.
\end{lem}
{\bf Proof.}  The loops discussed in Corollary \ref{sumup} generate a
subgroup of $G$ which factors as a product of subgroups of $G^{\rm
hyper}$ and $G^{g-g_0}$, since the loops preserve the decomposition
along $\partial U$ and act trivially on one side or the other. The
subgroup generated by these loops is also a subgroup of $G^{\rm
hyper}\times G^{g-g_0}$.  However, these loops generate inside
$G^{g-g_0}$ the diagonal copy of $\op{Sp}(2(g-g_0))$ acting on
$H^{1}(X_{o},{\mathbb C})^{g-g_0}\cong {\mathbb C} ^{2n^2(g-g_0)}$ (cf
the discussion above), 
and for the hyperelliptic case by the result of the previous section,
the group
\[
\op{Sp}(2g_0) \times \prod_{u \neq (\bone,0)}
\op{Sp}(H^{1}(X_{o},{\mathbb 
W}_{u})^{\rm hyper}) \subset G^{\rm hyper} .
\]
For general reasons, neither $G^{g-g_0}$ nor $G^{\rm hyper}$ can be
any bigger than these subgroups generated 
by our loops.  Therefore our loops generate all of $G^{g-g_0}$
(respectively $G^{\rm hyper}$), and we obtain that  
\[
G^{\rm hyper}\times G^{g-g_0} \subset G.
\]
The lemma is proven. \eop

\

\bigskip

\noindent
To finish the proof of Theorem~\ref{thm-lefschetz-pencil} we need some
elements of $G$ which mix up the factors
$H^{1}(X_{o},\op{ad}(\rho))^{\rm hyper}$ and
$H^{1}(X_{o},\op{ad}(\rho))^{g-g_0}$. We get these by the following an
``interchange of singularities'' argument, which basically comes down
to saying that certain singularities in the hyperelliptic family and
singularities in the complementary family are conjugate under the
global monodromy.

\begin{lem}
\label{interchange}
Let $\xi$ in $G^{\rm hyper}$ denote the monodromy element which acts
by the Dehn twist $D_{m-2,m-1}$ on the cohomology of the hyperelliptic
piece.  Let $\eta \in G^{g-g_0}$ denote a Dehn twist coming from a
double point on the complement of $U$. Then there is an element $\psi$
of the global monodromy group $G$ such that $\psi ^{-1}\xi \psi =
\eta$.
\end{lem}
{\bf Proof.} This is because the Dehn twists generating the
fundamental group of ${\mathbb P} H^0(Z,{\mathcal L} )-\bbD$ are
all conjugate since the discriminant divisor $\bbD$ is connected.
However, one must be a bit careful since we are looking at cohomology
with coefficients in a local system that is not necessarily preserved
by the full fundamental group: our representation $\rho$ is preserved
by a subgroup of finite index which corresponds to a covering of
${\mathbb P} H^0(Z,{\mathcal L} )$ ramified along $\bbD$, and in
this covering the inverse image of the discriminant locus might no
longer be connected.  We remedy this by stating somewhat more
explicitly how to construct $\psi$, but without actually writing down
the equations for this loop in ${\mathbb P} H^0(Z,{\mathcal L} )$
since that would be tedious.  Recall that we have assumed that the
last two branch points of the hyperelliptic curve were in branch cuts
where the passing transformation was the identity.  The curve acquires
a node when these two branch points come together; this node
corresponds to the monodromy element $\xi$.  Then this node can move
out of our coordinate patch $U$ and into the complementary region
$Z-U$.  At this point we are left with a hyperelliptic handle whose
equation is a small deformation of $x^{m-2}-y^2$ rather than
$x^m-y^2$.  The connectedness of the discriminant locus analogous to
$\bbD_\bbG-\bbG_{\infty}$ but for $m-2$ rather than $m$, and also for two
nodes at once, allows us to choose a path whereby the node which came
out of $U$ gets interchanged with another node corresponding to
$\eta$.  This connectedness statement, which holds for $k$ large
enough, is the second statement referred to in condition (d) at the
start of the argument.  After interchanging the two singularities, go
backwards along the path to send the $\xi$ node back into $U$.  All of
this corresponds to a path in the subvariety of $\bbD$ corresponding
to sections with two nodes.  (Geometrically this subvariety is the
codimension two nodal locus of the discriminant variety; we choose a
path which interchanges the two sheets of the discriminant which come
together along the nodal locus.) Choose a path $\psi$ which is near to
this path, but in the complement of the discriminant.  This has the
effect of giving the conjugation above.  Finally, notice that all of
this took place in a region where the representation $\rho$ is
trivial, so we can follow $\rho$ along the path $\psi$, i.e. the path
$\psi$ lifts to a path in the ramified covering on which $\rho$ is
defined. \eop

\

\bigskip

\noindent
Let $\xi \in G^{\rm hyper}$, $\eta \in G^{g-g_0}$, and $\psi \in G$ be
the elements from the above lemma. Note that $\eta$ is a generating
symplectic transvection in $\op{Sp}(2(g-g_0))$. Of course $\eta$ no
longer acts as a symplectic transvection but as a direct sum of $n^2$
copies of a transvection on $H^{1}(X_{o},{\mathbb C})^{g-g_0}$.

Now the proof of Theorem~\ref{thm-lefschetz-pencil} will be a
consequence of the following statement.

\begin{lem}
\label{statement}
The global monodromy group is the full product 
\[
G = \op{Sp}(H^{1}(X_{o},{\mathbb C})) \times \prod_{u \neq (\bone,0)}
\op{Sp}(H^{1}(X_{o},{\mathbb W}_{u}), 
\]
acting on $H^{1}(X_{o},{\mathbb C})\oplus(\oplus_{u \neq (\bone,0)}
H^{1}(X_{o},{\mathbb W}_{u}))$ by the sum of the standard representations. 
\end{lem}
{\bf Proof.}  Note first that $G$ is contained in the product.  Its
image in the first factor is the full group $\op{Sp}(2g) =
\op{Sp}(H^{1}(X_{o},{\mathbb C}))$, which 
is just the usual statement (the Deligne-Kazhdan-Margulis theorem
again) for cohomology with the trivial coefficient system ${\mathbb C}$.

Look at one of the pieces $H^{1}(X_{o},{\mathbb W}_{u})$, proceeding
in the spirit of the hyperelliptic discussion in the previous section
(the proof of Lemma~\ref{lem-irreducible} and
Theorem~\ref{thm-finite-index}). First we show that the action of $G$
on $H^{1}(X_{o},{\mathbb W}_{u})$ is irreducible.  For this, it
suffices to consider the action of the group algebra.  There is a
vector $v$ in $H^{1}(X_{o},{\mathbb W}_{u})^{\rm hyper}$ such that
\[
{\mathbb C} [G^{\rm hyper}]\cdot v = H^{1}(X_{o},{\mathbb W}_{u})^{\rm hyper}.
\]
We claim that 
\[
{\mathbb C} [G]\cdot v = H^{1}(X_{o},{\mathbb W}_{u}).
\]
Let $A$ be the image of the element $\xi -1$. It is a two-dimensional
subspace of $H^{1}(X_{o},{\mathbb W}_{u})$, and in particular it is
contained in ${\mathbb C} [G]\cdot v$.  Thus there are elements $v_1$
and $v_2$ of ${\mathbb C} [G]\cdot v$ such that $(\xi -1)v_1$ and
$(\xi -1)v_2$ span $A$.  On the other hand, the image $B$ of
$\eta-1=\psi ^{-1} (\xi -1)\psi$ is an isomorphic image of $A$, a
two-dimensional space contained in $H^{1}(X_{o},{\mathbb W}_{u})^{g -
g_{0}}$. The isomorphism is 
given by
\[
\psi^{-1} : A \stackrel{\cong}{\rightarrow} B.
\]
The vectors 
\[
\psi ^{-1} (\xi -1)\psi (\psi ^{-1}v_1) \;\;
\mbox{and}\;\;
\psi ^{-1} (\xi -1)\psi (\psi ^{-1}v_2)
\]
span $B$. In particular $B$ is contained in ${\mathbb C} [G]\cdot v$.  

Now we can write $H^{1}(X_{o},{\mathbb W}_{u})^{g - g_{0}} = {\mathbb
C} ^{2(g-g_0)} \oplus {\mathbb C} ^{2(g-g_0)}$ and the image of the
action of $G$ contains the diagonal copy of $\op{Sp}(2(g-g_0))$ (this is
the image of $G^{g-g_0}$).  The subspace $B$ is transverse to the
above decomposition, in other words it contains one basis element in
each piece.  It follows that the translates of $B$ by elements of
$G^{g-g_0}$ span $H^{1}(X_{o},{\mathbb W}_{u})^{g - g_{0}}$.  Therefore
\[
H^{1}(X_{o},{\mathbb W}_{u})^{g - g_{0}}\subset {\mathbb C} [G]\cdot v,
\]
so putting this together with the above we get that ${\mathbb C}
[G]\cdot v = H^{1}(X_{o},{\mathbb W}_{u})$, so the action of $G$ on
$H^{1}(X_{o},{\mathbb W}_{u})$ is irreducible as claimed (to get this
last deduction we use the standard fact that $G$ is semisimple so its
action decomposes as a direct sum of irreducible pieces).

Next note that the image of $G$ acting on $H^{1}(X_{o},{\mathbb
W}_{u})$ is a simple group. This is because the group $G^{\rm hyper}$
acting on $H^{1}(X_{o},{\mathbb W}_{u})$ contains an element whose
Jordan normal form has a single Jordan block of length one (see the
argument of the previous section for the hyperelliptic case).  As
before this implies that the image is simple.

Next we show that the image of $G$ acting on $H^{1}(X_{o},{\mathbb
W}_{u})$ is the full symplectic group $\op{Sp}(H^{1}(X_{o},{\mathbb
W}_{u}))$.  This again is by the same argument as in the previous
section, using the fact that the image of $G$ acting on $W_u$ contains
a copy of $\op{Sp}(2(g-g_0))$, and noting that we can insure that $g-g_0$
is large enough (by choosing $k$ big).

Finally, complete the proof of Lemma \ref{statement} by noting that
inside $G^{\rm hyper}$ we can find an element which acts with
different eigenvalues on each of the different pieces
$H^{1}(X_{o},{\mathbb W}_{u})$, the same element as exhibited in the
previous section. As then, this implies that
\[
G = \op{Sp}(2g) \times \prod_{u \neq (\bone,0)} \op{Sp}(H^{1}(X_{o},{\mathbb
W}_{u})).
\]
This completes the proof of the lemma. \eop

\

\bigskip

\noindent
{\bf Proof of Theorem~\ref{thm-lefschetz-pencil}:} The dihedral
Heisenberg representation $\rho$ which we chose here corresponds to a
smooth point in $M_B(X_{o},n)$, fixed under a finite index subgroup of
the monodromy group.  The action of the monodromy group which fixes
$\rho$ on the tangent space $T_{[\rho]}M_B(X_{o},n)$ at $\rho$ is
exactly the action on $H^{1}(X_{o},\op{ad}(\rho))$ we considered
above. The Zariski closure $G$ as described in Lemma \ref{statement},
acts with an open orbit. Either from the general consideration in
Lemma~\ref{lem-condition-ii} or else just by inspection, this
monodromy group has no characters.  Therefore by
Corollary~\ref{cor-localization} we get that the monodromy action on
$M_B(X_{o},n)$ is Zariski dense. \eop

\section{Further remarks} \label{sec-further}

\subsection{Topologically irreducible families} \label{sec-irreducible}

The requirement that our families have relatively large monodromy was
used in an essential way in the proofs of
Theorems~\ref{thm-finite-index} and
\ref{thm-lefschetz-pencil}. However this requirement 
 seems to be more an artifact of the method of proof rather than a
real condition on the family $f : X \to B$ which is necessary for the
density of the monodromy action. In this section we briefly examine
some consequences of the density, which will allow us to probe the
necessity of the `large monodromy' condition.

Recall that a smooth family of curves $f : X \to B$ is called {\em
topologically irreducible} if and only if there is no finite
collection of disjoint embedded circles in $X_{o}$ which is preserved
by the geometric monodromy. We have the following simple

\begin{lem} \label{lem-top-irreducible} Let $f : X \to B$ be a family of
smooth curves, such that the monodromy action of $\pi_{1}(B,o)$ has a
Zariski dense orbit on the Betti moduli space $M_{B}(X_{o},n)$ for
some $n \geq 1$. Then $f : X \to B$ is topologically irreducible.
\end{lem}
{\bf Proof.} Assume that one can find simple disjoint loops $a_{1},
\ldots, a_{k} \subset X_{o}$ such that the collection $\{ a_{1},
\ldots, a_{k} \}$ of free homotopy classes on $X_{o}$ is preserved by
$\op{mon}(\pi_{1}(B,o)) \subset \op{Map}(X_{o})$. Then for every $N
\in {\mathbb Z}$ we have a well defined
$\op{mon}(\pi_{1}(B,o))$-invariant regular function $\psi_{N} : M_{B}(X_{o},n)
\to {\mathbb C}$ on $M_{B}(X_{o},n)$, given by $\psi_{N}([\rho]) := 
\op{Tr}(\prod_{i=1}^{k} (\rho(a_{i}))^{N})$. But clearly for some
$N$ the function $\psi_{N}$ will be non-constant and so
$\op{mon}(\pi_{1}(B,o))$ can not have a Zariski dense orbit on
$M_{B}(X_{o},n)$. The lemma is proven. \eop

\

\bigskip

\

In particular, all families of curves satisfying the hypothesis of
Theorem~~\ref{thm-finite-index} and \ref{thm-lefschetz-pencil} 
will be topologically irreducible. 
Recently C.~McMullen has shown that topological irreducibility holds
very generally:  
every non-isotrivial
holomorphic  family  of curves is topologically irreducible
\cite[Proof of Theorem~3.1]{mcmullen}. In particular, this corollary of
our main results was already known.

When we
started the current project we were hoping that topological
irreducibility will allow one to distinguish symplectic Lefschetz
pencils (whose topology tends to be much softer)  
from projective Lefschetz pencils. 
In the meantime however, Ivan Smith
succeeded in showing \cite{smith} that all symplectic Lefschetz
fibrations over ${\mathbb P}^{1}$ are topologically irreducible. We
still expect that the stronger property {\bf GZD} (or the open orbit
property from Theorem~\ref{thm-lefschetz-pencil}) will allow one to
distinguish projective from symplectic Lefschetz pencils. We hope to
return to examples of this type in a future paper.

McMullen's result (with the alternative symplectic proof by Smith when
the base is ${\mathbb P}^1$) is the only evidence we have for the
following conjectural generalization of
Theorems~\ref{thm-finite-index} and \ref{thm-lefschetz-pencil}:

\begin{con} \label{con-isotrivial} 
For a family $f : X \to B$ assume that
$\op{mon}(\pi_{1}(B,o))$ is not a finite group. Then:
\begin{description}
\item[{\em (i)}] There is no meromorphic function on
$M_{B}(X_{o},n)^{\op{an}}$ which is invariant under 
the action of $\monb{n}(\pi_{1}(B,o))$ 
(equivalently there is no meromorphic function on 
$M_{DR}(X_{o},n)^{\op{an}}$ which is 
$\mondr{n}(\pi_{1}(B,o))$-invariant);
\item[{\em (ii)}] In the case of $M_{B}(X_{o},n)$, considered with its
natural structure of an affine algebraic variety,  there exist a point $x_{B}
\in M_{B}(X_{o},n)$ so that the orbit
\[
\monb{n}(\pi_{1}(B,o))\cdot x_{B}  \subset M_{B}(X_{o},n)
\] 
is Zariski dense in $M_{B}(X_{o},n)$.
\end{description}
\end{con}

\

\noindent
Procesi's theorem \cite{procesi}
implies that the field of rational functions on $M_{B}(X_{o},n)$ is
generated by traces of evaluation maps for conjugacy classes of simple
loops on $X_{o}$. One might hope (although we didn't find an argument) 
that the field of
$\pi_{1}(B,o)$-invariant rational functions on $M_{B}(X_{o},n)$ is
similarly generated by the traces of evaluation maps at finite invariant
collections of simple loops. If this were the case then McMullen's theorem
would imply the validity of the variant of Conjecture~\ref{con-isotrivial}
concerning algebraic meromorphic functions. The property {\bf AGZD1} (i.e. the 
conjecture as it is stated using analytic meromorphic functions) 
would seem to remain more elusive.

\subsection{Points $\Gamma$-near to a finite representation}
\label{subsec-gamma-near}

We briefly describe here another variation on the basic
result. Essentially, we have constructed a finite-image
representation, the dihedral Schr\"{o}dinger representation $\rho$,
which corresponds to a smooth point in $M_B(X_o,n)$ and which turns
out to be sufficient in order to get the Zariski-denseness
property. In an attempt to better understand what is going on, we can
explore a bit further the sense in which $\rho$ is near the rest of
$M_{B}(X_o,n)$.

Suppose a finitely presented group $\Gamma$ acts on an affine variety
$M$, and suppose $p\in M$ is a closed point in the smooth set of $M$,
fixed by the action. We say that another point $q\in M$ is {\em
$\Gamma$-near to $p$} if $p$ lies in the closure of the orbit
$\overline{\Gamma \cdot q}$. Let $N=\text{\bf Near}(M,p,\Gamma
)\subset M$ denote the subset of points $q$ which are $\Gamma$-near to
$p$. It is $\Gamma$-invariant.  Let $T_pN$ denote its tangent cone at
$p$, defined to be the set of limits of secants to $M$ going from $p$
to points $q\in N$ which approach $p$ (the limits of secants may be
taken in any real embedding of $M$). Note that $T_pN\subset T_pM$ is
an invariant subset of the tangent space to $M$ at $p$.  An easy
argument similar to that of Corollary \ref{cor-localization} shows
that if $T_pN$ is Zariski-dense in $T_pM$ then $N$ is Zariski-dense in
$M$.

Suppose $\gamma \in \Gamma$. Let $R\subset T_pM$ denote the span of the
eigenvectors of $\gamma$ whose eigenvalues have norm $<1$. We claim that
$R\subset T_pN$. Choose a smooth submanifold $V\subset M$ tangent to $R$.
Let $D$ be a ball neighborhood of $p$ in $M$. Then the collection 
\[
\{ V^k:= \gamma ^{-k} (W\cap \gamma ^kD)\}
\]
is a collection of manifolds with boundaries lying in the boundary of
$D$, which are tangent to $R$ at the origin, and whose curvature is
bounded (the $V\cap \gamma ^kD$ all lie in a sector preserved by
$\gamma ^{-1}$ and in which $\gamma ^{-1}$ smooths things out).  Thus
these converge to a manifold $V^{\infty}$ which is preserved by the
action of $\gamma$ and on which $\gamma$ acts with all eigenvalues
$<1$. In particular, $V^{\infty}\subset N$ which shows that $R\subset
T_pN$.

\begin{lem} \label{lem-nearness}
Suppose $\rho \in M_B(X_o, n)$ is the dihedral Schr\"{o}dinger
representation we have considered above. Suppose that a group $\Gamma
= \pi _1(B,o)$ acts, satisfying one of the hypotheses of Theorem A or
Theorem B. Let $N=\text{\bf Near}(M_B(X_o, n), \rho , \Gamma )$.  Then
$T_{\rho}N$ is Zariski-dense in $T_{\rho}M_B(X_o,n)$.
\end{lem}
{\bf Proof.} Recall that
\[
T_{[\rho]}M_B(X_o,n) = H^{1}(X_{o},\op{ad}(\rho))
\]
and that we have a decomposition
\[
H^{1}(X_{o},\op{ad}(\rho)) = \left(
\bigoplus_{u }H^{1}(X_{o},{\mathbb W}_{u}) \right)
\]
such that the monodromy group $\Gamma$ is Zariski-dense in the product
$G = \prod_{u } G_u$
with $G_u=\op{Sp}(H^{1}(X_{o},{\mathbb
W}_{u}))$.
Each component of the decomposition corresponds to a weight one
variation
of Hodge structure over $B$ which is irreducible and not unitary
(since both Hodge subspaces are nontrivial).
Therefore $\Gamma$ actually lies in a real form which decomposes
\[
\Gamma \subset G _{\mathbb R}= \prod_{u } G_{u,\mathbb R}
\]
and the $G_{u,\mathbb R}$ are noncompact real forms of the symplectic
groups. (One cannot have a real component whose complexification
splits into two components, because that would be a
complex group considered as a real group, which is never of Hodge type.)
The real Zariski closure of $\Gamma$, i.e. the
intersection of all real algebraic subsets of $G$ containing $\Gamma$,
is $G_{\mathbb R}$, since anything smaller would lead to a smaller
complex Zariski closure.

Let $pr_u:G_{\mathbb R}\rightarrow G_{u,\mathbb R}$
denote the projection. Let $E_u\subset G_{\mathbb R}$ be the
real algebraic subset of elements $g$ such that all of the eigenvalues of
$pr_u(g)\in  G_{u,\mathbb R}$ have norm $1$.  This
is a proper subset since $G_{u,\mathbb R}$ is noncompact.  The union
$\bigcup _u E_u$ is again a proper real algebraic subset of $G_{\mathbb R}$,
so it cannot contain $\Gamma$.

Thus there is an element  $\gamma\in \Gamma$
such that every projection $pr_u(\gamma )$ has at least one eigenvalue of
norm different
from $1$. On the other hand these projections have determinant one, so each
$pr_u(\gamma )$ has at least one eigenvalue of norm $<1$.

In particular, there is a vector $v\in
H^{1}(X_{o},\op{ad}(\rho))$ such that $v$ is in the span of the
eigenvectors of $\gamma$ corresponding to eigenvalues of absolute
value $<1$, and such that $v$ has a nonzero component in all of the
irreducible factors of $H^{1}(X_{o},\op{ad}(\rho))$. Using the fact
that the complex
Zariski closure of $\Gamma$ contains a product of symplectic groups
(corresponding to the decomposition into irreducible pieces of the
representation $H^{1}(X_{o},\op{ad}(\rho))$), we find that the orbit of
the vector $v$ under the action of $\Gamma$ is Zariski-dense in
$T_{[\rho]}M_B(X_o,n)$.  On the other hand, from our discussion prior
to the present lemma, $v$ is in the $\Gamma$-invariant subset
$T_{\rho}N$.  Thus $T_{\rho}N$ is Zariski-dense.  \eop

\begin{cor} \label{cor-nearness}
The moduli space $M_B(X_o, n)$ contains a smooth point, the dihedral
\linebreak
Schr\"{o}dinger representation $\rho$, such that the set of points
$\text{\bf Near}(M_B(X_o, n), \rho , \Gamma )\subset M_B(X_o, n)$ which 
are $\Gamma$-near  to
$\rho$, is Zariski-dense in $M_B(X_o,n)$.
\end{cor}

It is clear that any $\Gamma$-invariant regular (i.e. holomorphic algebraic)
function takes the 
same value at $\rho$ as at every point of $\text{\bf Near}(M_B(X_o, n), 
\rho , \Gamma )$. 
In particular, this corollary implies the result that there are no
$\Gamma$-invariant regular functions. This result is weaker than our main results
about nonexistence of invariant meromorphic functions, but does provide a slightly
different conceptual route to the
topological irreducibility result refered to above. 

Our motivation for introducing the notion of $\Gamma$-nearness is that
the first and second authors asked some time ago whether there was any
sense in which the finite-image representations could take up a big
place in $M_B(X_o,n)$. The short answer to that question is that, by
Jordan's theorem, the finite image representations occupy a rather
small place in that there are only finitely many outside of
representations which factor through a normalizer of a torus (and
those which factor in this way lie in a closed subset of relatively
high codimension).  However, Corollary \ref{cor-nearness} provides the
slightly more subtle answer that, in the presence of a large monodromy
action, if you start out very near to a certain finite-image
represenation (such as one of our dihedral Schr\"{o}dinger
representations) and then let the monodromy act, then you can get out
to a significant part of $M_B(X_o,n)$.

\subsection{Other groups}
\label{subsec-othergroups}

Finally, let us explicitly state that we expect that all of the
results and conjectures of this paper to hold for Betti and deRham
cohomology with coefficients in an arbitrary complex reductive group
$G$. Specifically we make the following

\begin{con} \label{con-G} Let $f : X \to B$ be a smooth algebraic
family of curves and $G$ a complex reductive group. Then:
\begin{description}
\item[{\em (i)}] Assume that  $f$ is not isotrivial. 
Then the families 
\[
M_{B}(X/B,G) \to B \quad \text{and} \quad   M_{DR}(X/B,G)
\to B
\]
of relative Betti and de Rham cohomology with coefficients in
$G$ are {\bf GZD} when  equipped with the non-abelian Gauss-Manin
connection.
\item[{\em (ii)}] Assume that $f$ comes from a projective Lefschetz
pencil of sufficiently  high degree. Then there exists a smooth point
$\rho \in M_{B}(X_{o},G)$ which is fixed by a finite index subgroup
$\Gamma \subset \pi_{1}(B,o)$ and for which the Zariski closure of
$\Gamma$ in $GL(T_{\rho}M_{B}(X_{o},G))$ acts with an open orbit on 
$T_{\rho}M_{B}(X_{o},G)$.
\end{description}
\end{con}


\begin{thebibliography}{}


\bibitem[\protect\citeauthoryear{Deligne}{72}]{deligne-hodge2}
P.~Deligne.
\newblock {Th\'{e}orie} de {Hodge} {II}.
\newblock {\em Publications Math\'{e}matiques de l'{I.}{H.}{E.}{S.}}, 40:5--57,
  1972.

\bibitem[\protect\citeauthoryear{Deligne}{73}]{deligne-ode-book}
P.~Deligne.
\newblock {\em \'{E}quations diff\'{e}rentielles \`{a} points singuliers
  r\'{e}guliers}.
\newblock Springer-Verlag, 1973.
\newblock Lecture Notes in Mathematics, Vol 163.

\bibitem[\protect\citeauthoryear{Deligne}{80}]{deligne-weil2}
P.~Deligne.
\newblock Le conjecture de {W}eil {I}{I}.
\newblock {\em Publ. Math. {I.}{H.}{E.}{S.}}, 52:313--428, 1980.

\bibitem[\protect\citeauthoryear{Gallo \bgroup \em et al.\egroup }{00}]{gkm}
D.~Gallo, M.~Kapovich, and A.~Marden.
\newblock The monodromy groups of {S}chwarzian equations on closed {R}iemann
  surfaces.
\newblock {\em Ann. of Math. (2)}, 151(2):625--704, 2000.

\bibitem[\protect\citeauthoryear{Goldman}{97}]{goldman-ergodic}
W.~Goldman.
\newblock Ergodic theory on moduli spaces.
\newblock {\em Ann. of Math. (2)}, 146(3):475--507, 1997.

\bibitem[\protect\citeauthoryear{Grothendieck}{68}]{grothendieck-crystals}
A.~Grothendieck.
\newblock Crystals and the de {Rham} cohomology of schemes.
\newblock In {\em Dix expos\'{e}s sur la cohomologie des sch\'{e}mas}.
  North-Holland, Amsterdam, 1968.

\bibitem[\protect\citeauthoryear{Janssen}{83}]{janssen}
W.A.M. Janssen.
\newblock Skew-symmetric vanishing lattices and their monodromy groups.
\newblock {\em Math. Ann.}, 266(1):115--133, 1983.

\bibitem[\protect\citeauthoryear{Katz}{68}]{katz-sga} 
N.~Katz
\newblock Etude cohomologique de pinceaux de Lefschetz.  
\newblock Expose XVIII in {\em S\'{e}minaire de G\'{e}om\'{e}trie
Alg\'{e}brique 
du Bois-Marie 1967--1969 (SGA 7 II).} Dirig\'{e} par P.Deligne et
N.Katz. Lecture Notes in Mathematics, Vol. 340.


\bibitem[\protect\citeauthoryear{Looijenga}{96}]{looijenga-notes}
E.~Looijenga.
\newblock Cohomology and intersection homology of algebraic varieties.
\newblock In {\em Complex algebraic geometry (Park City, UT, 1993)}, pages
  221--263. Amer. Math. Soc., Providence, RI, 1996.

\bibitem[\protect\citeauthoryear{Looijenga}{97}]{looijenga}
E.~Looijenga.
\newblock Prym representations of mapping class groups.
\newblock {\em Geom. Dedicata}, 64(1):69--83, 1997.

\bibitem[\protect\citeauthoryear{Lubotzky-Magid}{85}]{repsfg}
A.~Lubotzky and A.~Magid.
\newblock Varieties of representations of finitely generated groups.
\newblock {\em Mem. Amer. Math. Soc.}, 58(336):xi+117, 1985.


\bibitem[\protect\citeauthoryear{McMullen}{00}]{mcmullen}
C.~McMullen.
\newblock From dynamics on surfaces to rational points on curves.
\newblock {\em Bull. Amer. Math. Soc. (N.S.)} 37:119--140, 2000.

\bibitem[\protect\citeauthoryear{Mess}{92}]{mess}
G.~Mess.
\newblock A note on hyperelliptic curves.
\newblock {\em Proc. Amer. Math. Soc.}, 115:849-852, 1992.


\bibitem[\protect\citeauthoryear{Mumford}{93}]{mumford} D.~Mumford
\newblock Tata lectures on theta. III.  
\newblock With the
collaboration of Madhav Nori and Peter Norman. Progress in
Mathematics, 97.  Birkh\"{a}user Boston, Inc., Boston, MA, 1991.

\bibitem[\protect\citeauthoryear{Procesi}{74}]{procesi}
C.~Procesi.
\newblock Finite dimensional representations of algebras.
\newblock {\em Israel J. Math.}, 19:169--182, 1974.

\bibitem[\protect\citeauthoryear{Smith}{01}]{smith}
I.~Smith.
\newblock {Geometric monodromy and the hyperbolic disc}, 
\newblock {\em Q. J. Math.} 52:217--228, 2001.

\bibitem[\protect\citeauthoryear{Simpson}{95}]{simpson-moduli2}
C.~Simpson.
\newblock Moduli of representations of the fundamental group of a smooth
  projective variety - {I}{I}.
\newblock {\em Publications Math\'{e}matiques de l' {I.}{H.}{E.}{S.}},
  80:5--79, 1995.

\bibitem[\protect\citeauthoryear{Tibar}{02}{a}]{Tibar1}
M.~Tibar. 
\newblock Connectivity via nongeneric pencils.
\newblock {\em Internat. J. Math.}{\bf 13} (2002), 111-123.

\bibitem[\protect\citeauthoryear{Tibar}{02}{b}]{Tibar2}
M.~Tibar. 
\newblock Vanishing cycles of pencils of hypersurfaces. 
\newblock Preprint {\tt math.AG/0204094}.





\end{thebibliography}
\end{document}